\documentclass[fleqn]{article}

\usepackage{amscd}
\usepackage{amsthm}
\usepackage{amsmath}
\usepackage{amssymb}
\usepackage{xypic}
\usepackage[latin5]{inputenc}
\usepackage[english]{babel}

\newtheorem{thm}{Theorem}[section]
\newtheorem{cor}[thm]{Corollary}
\newtheorem{lem}[thm]{Lemma}
\newtheorem{prop}[thm]{Proposition}
\theoremstyle{definition}
\newtheorem{defn}[thm]{Definition}
\newtheorem{example}[thm]{Example}
\theoremstyle{remark}
\newtheorem{rem}[thm]{Remark}

\numberwithin{equation}{subsection}

\begin{document}

\newcommand{\Dif}{\protect \mbox{Diff}\,}

\newcommand{\dif}{\ensuremath{{\protect \mbox{Diff}_{\scriptsize x_{\tiny
0}}}}}

\newcommand{\difk}{\ensuremath{{\protect \mbox{Diff}^k_{\scriptsize x_{\tiny
0}}}}}
\newcommand{\difrk}{\ensuremath{{\protect \mbox{Diff}^{r+k}_{\scriptsize x_{\tiny
0}}}}}
\newcommand{\difrmasuno}{\ensuremath{{\protect \mbox{Diff}^{r+1}_{\scriptsize x_{\tiny
0}}}}}

\newcommand{\Hom}[3]{\ensuremath{{\protect \mbox{Hom}_{\scriptsize
#1}(#2,#3)}}}
\newcommand{\fib}[1]{\ensuremath{#1}_{\scriptsize x_{\tiny 0}}}
\newcommand{\kom}[1]{\ensuremath H_{\scriptsize x_{\tiny 0}}^{#1}}
\newcommand{\moduli}{\ensuremath{{\protect \mathfrak{M}_n^r}}}

\title{\textsc{Moduli spaces for finite-order jets of Riemannian
metrics}}
\author{Gordillo, A., Navarro, J. and Sancho, J.B.}
\date{\today}

\maketitle

\begin{abstract}
We construct the moduli space of $\,r-$jets of Riemannian metrics
at a point on a smooth manifold. The construction is closely
related to the problem of classification of jet metrics via
differential invariants.

The moduli space is proved to be a differentiable space which
admits a finite canonical stratification into smooth manifolds. A
complete study on the stratification of moduli spaces is carried
out for metrics in dimension $\,n=2\,$.
\end{abstract}

\bigskip

\section*{Introduction}

Let $\,X\,$ be an $\,n-$dimensional smooth manifold. Fixed a point
$\,x_0\in X\,$ and an integer $\,r\geq0\,$, we will denote by
$\,J_{x_0}^rM\,$ the smooth manifold of $\,r-$jets at $\,x_0\,$ of
Riemannian metrics on $\,X\,$. On the manifold $\,J_{x_0}^rM\,$,
there exists a natural action of the group $\,\mbox{Diff}_{x_0}\,$
of germs at $\,x_0\,$ of local diffeomorphisms leaving $\,x_0\,$
fixed, so it yields an equivalence relation on $\,J_{x_0}^rM\,$:
$$
j_{x_0}^rg\,\equiv\,j_{x_0}^r\bar{g} \,\,\,\Longleftrightarrow
\,\,\,
j_{x_0}^r(\tau^*g)=j_{x_0}^r\bar{g}\,\,,\,\,\mathrm{for\,\,some\,\,}\tau\in
\mathrm{Diff}_{x_0}\,.
$$

The quotient space
$\,\mathbb{M}_n^r:=J_{x_0}^rM/\mbox{Diff}_{x_0}\,$ is called
moduli space for $\,r-$jets of Riemannian metrics in dimension
$\,n\,$. It depends neither on the point $\,x_0\,$ nor on the
$n-$dimensional manifold $\,X\,$ chosen.

The purpose of this paper is to study the structure of moduli
spaces $\,\mathbb{M}_n^r\,$.

\medskip

Moduli spaces $\,\mathbb{M}_n^r\,$ have been studied in the
literature through their function algebras
$\,\mathcal{C}^{\infty}(\mathbb{M}_n^r):=\mathcal{C}^{\infty}(J_{x_0}^rM)^{\textrm{Diff}_{x_0}}\,$.
This function algebra $\,\mathcal{C}^{\infty}(\mathbb{M}_n^r)\,$
is nothing but the algebra of differential invariants of order
$\,\leq r\,$ of Riemannian metrics. Mu{\~n}oz and Vald\'es
(\cite{masval},\cite{masval2}) prove that it is an essentially
finitely-generated algebra and they determine the number of its
functionally independent generators. In a more general setting,
Vinogradov (\cite{vinogradov}) has pointed out a simple and
natural relationship between the algebra of differential
invariants of homogeneous geometric structures and their
characteristic classes. (See also \cite{vinogradovII}.)

Let us also mention that in \cite{garmas} Garc\'{\i}a and Mu{\~n}oz
obtain a moduli space for linear frames, which has structure of
smooth manifold.

However, apart from some trivial exceptions, moduli spaces
$\,\mathbb{M}_n^r\,$ of jet metrics are not smooth manifolds, but
they possess a differentiable structure in a more general sense:
that of a differentiable space. (The typical example of
differentiable space is a closed subset $\,Y\subseteq
\mathbb{R}^m\,$ where a function $\,f:Y\rightarrow \mathbb{R}\,$
is said to be differentiable if it is the restriction to $\,Y\,$
of a smooth function on $\,\mathbb{R}^m\,$, see \cite{navsan}.)

In addition, the differentiable structure of $\,\mathbb{M}_n^r\,$
is not too far from a smooth structure, since it admits a
stratification by a finite number of smooth submanifolds. Our
results can be summed up in the following

\smallskip

\begin{thm}\label{teoremaintroduccion} Every moduli space $\,\mathbb{M}_n^r\,$ is a
differentiable space and it admits a finite canonical
stratification
$$
\mathbb{M}_n^r=S_{[H_0]}^r \sqcup \ldots \sqcup S_{[H_s]}^r\,,
$$
for locally closed subspaces $\,S_{[H_i]}^r\,$ which are smooth
manifolds. Moreover, one of them is an open connected dense subset
of $\,\mathbb{M}_n^r\,$.
\end{thm}

\medskip

Each stratum of this decomposition of the space
$\,\mathbb{M}_n^r\,$ consists of those jet metrics having
essentially the same group of automorphisms. To be more precise,
let us denote by $\,[H]\,$ the conjugacy class of a closed
subgroup $\,H\,$ of the orthogonal group $\,O(n)\,$. Then
$\,S_{[H]}^r\,$ is the set of equivalence classes of jet metrics
$\,j_{x_0}^rg\,$ whose group of automorphisms
$\,\mbox{Aut}(j_{x_0}^rg)\,$ is conjugate to $\,H\,$, viewing
$\,\mbox{Aut}(j_{x_0}^rg)\,$ as a subgroup of the orthogonal group
$\,O(T_{x_0}X,g_{x_0})\simeq O(n)\,$.

It is convenient to notice that Theorem \ref{teoremaintroduccion}
is not valid for semi-Riemannian metrics. For metrics of any
signature, the problem lies on the existence of non-closed orbits
for the action of $\,\mbox{Diff}_{x_0}\,$ on the space
$\,J_{x_0}^rM\,$ of $\,r-$jets of such metrics, which means that
the corresponding moduli space $\,J_{x_0}^rM/\mbox{Diff}_{x_0}\,$
is not a $\,T_1\,$ topological space, and consequently, it does
not admit a structure of differentiable space either.

In dimension $\,n=2\,$, we improve the above theorem by
determining exactly all the strata which appear in the
decomposition of each moduli space $\,\mathbb{M}_{n=2}^r\,$. Let
us consider the only, up to conjugacy, closed subgroups of the
orthogonal group $\,O(2)\,$: the finite group $\,K_m\,$ of
rotations of order $\,m\,$ ($\,m\geq 1\,$), the dihedral group
$\,D_m\,$ of order $\,2m\,$ ($\,m\geq 1\,$), the special
orthogonal group $\,SO(2)\,$ and $\,O(2)\,$ itself. The
stratification of $\,\mathbb{M}_2^r\,$ is determined by the
following

\smallskip

\begin{thm} The strata in the moduli space
$\,\mathbb{M}_{n=2}^r\,$ correspond exactly to the following
conjugacy
classes:$\,[O(2)]\,,\,[D_1]\,,\ldots,\,[D_{r-2}]\,,\,[K_1]\,,\,\ldots\,,\,[K_{r-4}]\,$.
(And also $\,[K_1]\,$, if $\,r = 4\,$.)
\end{thm}

\medskip

Finally, we include two appendices. In the first one, we give a
brief discussion of the notion of differential invariant. In the
second one, we analyze the equivalence problem for infinite-order
jets of Riemannian metrics.

\bigskip

\section{Preliminaries}\label{preliminares}

\subsection{Quotient spaces}

Throughout this paper, we are going to handle geometric objects of
a more general nature than smooth manifolds, which appear when one
considers the quotient of a smooth manifold by the action of a Lie
group.

\begin{defn}
Let $\,X\,$ be a topological space. A \textbf{sheaf of continuous
functions} on $\,X\,$ is a map $\,\mathcal{O}_X\,$ which assigns a
subalgebra $\,\mathcal{O}_X(U)\subseteq
\mathcal{C}(U,\mathbb{R})\,$ to every open subset $\,U\subseteq
X\,$, with the following condition:

For every open subset $\,U\subseteq X\,$, every open cover
$\,U=\bigcup U_i\,$ and every function $\,f:U\rightarrow
\mathbb{R}\,$, it is verified
$$
f\in \mathcal{O}_X(U)\,\, \Longleftrightarrow \,\, f|_{U_i}\in
\mathcal{O}_X(U_i)\,,\,\,\,\forall\, i\,.
$$
\end{defn}

\medskip

In particular, if $\,V\subseteq U\,$ are open subsets in $\,X\,$,
then it is verified
$$
f\in \mathcal{O}_X(U)\,\, \Longrightarrow \,\, f|_V\in
\mathcal{O}_X(V)\,.
$$

\smallskip

\begin{defn}
We will call \textbf{ringed space} the pair
$\,(X,\mathcal{O}_X)\,$ formed by a topological space $\,X\,$ and
a sheaf of continuous functions $\,\mathcal{O}_X\,$ on $\,X\,$.
\end{defn}

\medskip

Although the concept of ringed space in the literature, specially
in that concerning Algebraic Geometry, is much broader, the
previous definition is good enough for our purposes.

\smallskip

Every open subset $\,U\,$ of a ringed space
$\,(X,\mathcal{O}_X)\,$ is itself, in a very natural way, a ringed
space, if we define $\,\mathcal{O}_U(V):=\mathcal{O}_X(V)\,$ for
every open subset $\,V\subseteq U\,$.

Hereinafter, a ringed space $\,(X,\mathcal{O}_X)\,$ will usually
be denoted just by $\,X\,$, dropping the sheaf of functions.

\smallskip

\begin{defn}
Given two ringed spaces $\,X\,$ and $\,Y\,$, a \textbf{morphism of
ringed spaces} $\,\varphi:X\rightarrow Y\,$ is a continuous map
such that, for every open subset $\,V\subseteq Y\,$, the following
condition is held:
$$
f\in \mathcal{O}_Y(V)\,\, \Longrightarrow \,\, f\circ\varphi \in
\mathcal{O}_X(\varphi^{-1}(V))\,.
$$
A morphism of ringed spaces $\,\varphi:X\rightarrow Y\,$ is said
to be an \textbf{isomorphism} if it has an inverse morphism, that
is, there exists a morphism of ringed spaces $\,\phi:Y\rightarrow
X\,$ verifying $\,\varphi \circ \phi=\mbox{Id}_Y\,$, $\,\phi \circ
\varphi=\mbox{Id}_X\,$.
\end{defn}

\medskip

\begin{example} \textbf{(Smooth manifolds)} The space
$\,\mathbb{R}^n\,$, endowed with the sheaf
$\,\mathcal{C}_{\mathbb{R}^n}^{\infty}\,$ of smooth functions, is
an example of ringed space. An $\,n-$smooth manifold is precisely
a ringed space in which every point has an open neighbourhood
isomorphic to
$\,(\mathbb{R}^n,\mathcal{C}_{\mathbb{R}^n}^{\infty})\,$. Smooth
maps between smooth manifolds are nothing but morphisms of ringed
spaces.
\end{example}

\smallskip

\begin{example} \textbf{(Quotients by the action of a Lie group)}
Let $\,G\times X\rightarrow X\,$ be a smooth action of a Lie group
$\,G\,$ on a smooth manifold $\,X\,$, and let $\,\pi:X\rightarrow
X/G\,$ be the canonical quotient map.

We will consider on the quotient topological space $\,X/G\,$ the
following sheaf $\,\mathcal{C}_{X/G}^{\infty}\,$ of
``differentiable'' functions:

For every open subset $\,V\subseteq X/G\,$,
$\,\mathcal{C}_{X/G}^{\infty}(V)\,$ is defined to be
$$
\mathcal{C}_{X/G}^{\infty}(V):=\{f:V\longrightarrow \mathbb{R}: f
\circ \pi \in \mathcal{C}^{\infty}(\pi^{-1}(V))\}\,.
$$
Note that there exists a canonical $\,\mathbb{R}-$algebra
isomorphism:
$$
\begin{CD}
\mathcal{C}_{X/G}^{\infty}(V) @=
\mathcal{C}^{\infty}(\pi^{-1}(V))^G\\
f & \longmapsto & f\circ \pi\,.
\end{CD}
$$

The pair $\,(X/G,\mathcal{C}_{X/G}^{\infty})\,$ is an example of
ringed space, which we will call \textbf{quotient ringed space} of
the action of $\,G\,$ on $\,X\,$.

As it would be expected, this space verifies the \textbf{universal
quotient property}: Every morphism of ringed spaces
$\,\varphi:X\rightarrow Y\,$, which is constant on every orbit of
the action of $\,G\,$ on $\,X\,$, factors uniquely through the
quotient map $\,\pi:X\rightarrow X/G\,$, that is, there exists a
unique morphism of ringed spaces $\,\tilde{\varphi}:X/G\rightarrow
Y\,$ verifying $\,\varphi=\tilde{\varphi}\circ \pi\,$.
\end{example}

\smallskip

\begin{example} \textbf{(Inverse limit of smooth manifolds)} Sometimes we will
consider an inverse system
$$
\cdots \longrightarrow X_{r+1} \longrightarrow X_r \longrightarrow
\cdots \longrightarrow X_1
$$
of smooth mappings between smooth manifolds (or, with some more
generality, an inverse system of ringed spaces).

The inverse limit $\,{\displaystyle \lim_{\leftarrow}}\, X_r\,$ is
a ringed space in the following natural way. On $\,{\displaystyle
\lim_{\leftarrow}}\, X_r\,$ it is considered the inverse limit
topology, that is, the initial topology induced by the evident
projections $\,p_s:\,{\displaystyle \lim_{\leftarrow}}\,
X_r\rightarrow X_s\,$. A real function on an open subset of
$\,{\displaystyle \lim_{\leftarrow}}\, X_r\,$ is said to be
``differentiable'' if it locally coincides with the composition of
a projection $\,p_s:{\displaystyle \lim_{\leftarrow}}\, X_r
\rightarrow X_s\,$ and a smooth function on $\,X_s\,$.

The topological space $\,{\displaystyle \lim_{\leftarrow}}\,
X_r\,$ endowed with the above sheaf of differentiable functions is
a ringed space satisfying the suitable universal property:

For every ringed space $\,Z\,$, there exists the bijection
$$
\begin{CD}
\mathrm{Hom}\,(Z,{\displaystyle \lim_{\leftarrow}}\, X_r) @=
{\displaystyle \lim_{\leftarrow}}\,
\mathrm{Hom}\,(Z,X_r)\\
\varphi & \longmapsto & (\ldots,p_r\circ \varphi,\ldots)\,.
\end{CD}
$$
\end{example}

\smallskip

\begin{example}
Let $\,Z\,$ be a locally closed subspace of $\,\mathbb{R}^n\,$. We
define the sheaf $\,\mathcal{C}_Z^{\infty}\,$ of differentiable
functions on $\,Z\,$ to be the sheaf of functions locally
coinciding with restrictions of smooth functions on
$\,\mathbb{R}^n\,$. The pair $\,(Z,\mathcal{C}_Z^{\infty})\,$ is
another example of ringed space.
\end{example}

\medskip

\begin{defn}
A (reduced) \textbf{differentiable space} is a ringed space in
which every point has an open neighbourhood isomorphic to a
certain locally closed subspace $\,(Z,\mathcal{C}_Z^{\infty})\,$
in some $\,\mathbb{R}^n\,$.

A map between differentiable spaces is called
\textbf{differentiable} if it is a morphism of ringed spaces.
\end{defn}

\smallskip

\begin{thm} {\rm \textbf{(Schwarz \cite{schwarz},\cite{navsan} Th.
11.14)}}\label{schwarz} Let $\,G\rightarrow \mbox{Gl}(V)\,$ be a
finite-dimensional linear representation of a compact Lie group
$\,G\,$. The quotient space $\,V/G\,$ is a differentiable space.
\end{thm}

\medskip

More precisely: Let $\,p_1,\ldots,p_s\,$ be a finite set of
generators for the $\,\mathbb{R}-$algebra of $\,G-$invariant
polynomials on $\,V\,$; these invariants define an isomorphism of
ringed spaces
$$
\xymatrix{(p_1,\ldots,p_s)\,:\,V/G \ar@{=}[r] & Z\subseteq
\mathbb{R}^s\,,}
$$
$Z\,$ being a closed subspace of $\,\mathbb{R}^s\,$.

\bigskip

\subsection{Normal tensors}

Let $\,X\,$ be an $\,n-$dimensional smooth manifold. Fix a point
$\,x_0 \in X\,$ and a semi-Riemannian metric $\,g\,$ on $\,X\,$ of
fixed signature $\,(p,q)\,$, with $\,n=p+q\,$. Let us recall
briefly some definitions and results:

\smallskip

\begin{defn} A coordinate system $\,(z_1, \ldots , z_n)\,$ in a
neighbourhood of $\,x_0\,$ is said to be a \textbf{normal
coordinate} system for $\,g\,$ at the point $\,x_0\,$ if the
geodesics passing through $\,x_0\,$ at $\,t=0\,$ are precisely the
``straight lines'' $\,\{ z_1(t) = \lambda_1 t , \ldots , z_n(t) =
\lambda_n t \}\,$, where $\,\lambda_i \in \mathbb{R}\,$.
\end{defn}

\medskip

In particular, $\,x_0\,$ is the origin of any normal coordinate
system for $\,g\,$ at $\,x_0\,$.

\smallskip

\begin{rem}
Observe that we do not require $\,(\partial_{z_1} , \ldots ,
\partial_{z_n})\,$ to be an orthonormal basis of $\,T_{x_0}X\,$.
\end{rem}

\medskip

As it is well known, via the exponential map
$\,\mbox{exp}_g:T_{x_0}X\rightarrow X\,$, normal coordinate
systems on $\,X\,$ correspond bijectively to linear coordinate
systems on $\,T_{x_0}X\,$. Therefore, two normal systems differ in
a linear coordinate transformation.

\smallskip

\begin{prop}\label{exponencial} Let $\,g\,$, $\,\bar{g}\,$ be two
semi-Riemannian metrics on $\,X\,$. Let us also consider their
corresponding exponential maps
$\,\mathrm{exp}_g\,,\,\mathrm{exp}_{\bar{g}}:T_{x_0}X\rightarrow
X\,$. For every $\,r\geq 0\,$ it is verified:
$$
j_{x_0}^rg=j_{x_0}^r\bar{g}\,\, \Longrightarrow \,\,
j_0^{r+1}(\mathrm{exp}_g)=j_0^{r+1}(\mathrm{exp}_{\bar{g}})\,.
$$
\end{prop}

\medskip

As a consequence of Proposition \ref{exponencial}, whose proof is
routine, normal coordinate systems at $\,x_0\,$ for a metric
$\,g\,$ are determined up to the order $\,r+1\,$ by the jet
$\,j_{x_0}^rg\,$. This fact will be used later on with no more
explicit mention.

\smallskip

\begin{defn} Let $\,r \geq 1\,$ be a fixed integer and let $\,x_0 \in X\,$. The space of
\textbf{normal tensors} of order $\,r\,$ at $\,x_0\,$, which we
will denote by $\,N_r\,$, is the vector space of
$\,(r+2)-$covariant tensors $\,T\,$ at $\,x_0\,$ having the
following symmetries:

- $\,T\,$ is symmetric in the first two and last $\,r\,$ indices:
$$T_{ijk_1 \ldots k_r} = T_{jik_1 \ldots k_r} \quad , \quad T_{ijk_1
\ldots k_r} = T_{ijk_{\sigma(1)} \ldots k_{\sigma(r)}}\quad, \quad
\forall \ \sigma \in S_r\,;$$

- the cyclic sum over the last $\,r+1\,$ indices is zero:
$$T_{ijk_1 \ldots k_r} + T_{ik_rjk_1 \ldots k_{r-1}} + \ldots +
T_{ik_{1} \ldots k_{r}j} = 0\,.$$

If $\,r=0\,$, we will assume $\,N_0\,$ to be the set of
semi-Riemannian metrics at $\,x_0\,$ of a fixed signature
$\,(p,q)\,$ (which is an open subset of $\,S^2T^*_{x_0}X\,$, but
not a vector subspace).
\end{defn}

\medskip

A simple computation shows that, in general, $\,N_1 = 0\,$.
Moreover, in \cite{eps} it is proved that $\,N_r\,$ ($r\geq 2$) is
a linear irreducible representation of the linear group
$\,\mbox{Gl}\,(T_{x_0}X)\,$.

To show how a semi-riemannian metric $\,g\,$ produces a sequence
of normal tensors $\,g^r_{x_0}\,$ at $\,x_0\,$, let us recall this
classical result:

\smallskip

\begin{lem}{\bf (Gauss Lemma)}\label{Gauss} Let $\,(z_1,\ldots,z_n)\,$ be germs of
coordinates centred at $\,x_0\in X\,$. These coordinates are
normal for the germ of a semi-Riemannian metric $\,g\,$ if and
only if the metric coefficients $\,g_{ij}\,$ verify the equations
$$
\sum_j g_{ij}z_j= \sum_j g_{ij}(x_0) z_j\,.
$$
\end{lem}

\medskip

Let $\,(z_1,\ldots , z_n)\,$ be a normal coordinate system for
$\,g\,$ at $\,x_0 \in X\,$ and let us denote:
$$g_{ij,k_1 \ldots k_r} := \frac{\partial^r g_{ij}}{\partial z_{k_1} \ldots \partial z_{k_r}}
(x_0)\,.$$

If we differentiate $\,r+1\,$ times the identity of the Gauss
Lemma, we obtain:
$$g_{ik_0,k_1 \ldots k_r}+g_{ik_1,k_2 \ldots k_r k_0} + \cdots + g_{ik_r,k_0 \ldots k_{r-1}}\,=\,0\,.$$
This property, together with the obvious fact that the
coefficients $\,g_{ij,k_1 \ldots k_r}\,$ are symmetric in the
first two and in the last $\,r\,$ indices, allows to prove that
the tensor
$$g^r_{x_0} := \sum_{ijk_1 \ldots k_r} g_{ij,k_1
\ldots k_r} \, \text{d}z_i \otimes \text{d}z_j \otimes
\text{d}z_{k_1} \otimes \ldots \otimes \text{d}z_{k_r}$$ is a
normal tensor of order $\,r\,$ at $\,x_0 \in X\,$. This
construction does not depend on the choice of the normal
coordinate system $\,(z_1,\ldots , z_n)\,$.

\smallskip

\begin{defn} The tensor $\,g^r_{x_0}\,$ is called the $\,r-$th \textbf{normal tensor of
the metric} $\,g\,$ at the point $\,x_0\,$.

As a consequence of $\,N_1 = 0\,$, the first normal tensor of a
metric $\,g\,$ is always zero, $\,g^1_{x_0} = 0\,$.
\end{defn}

\medskip

The normal tensors associated to a metric were first introduced by
Thomas \cite{tho}. The sequence
$\,\{g_{x_0},g_{x_0}^2,g_{x_0}^3,\ldots,g_{x_0}^r\}\,$ of normal
tensors of the metric $\,g\,$ at a point $\,x_0\,$ totally
determines the sequence
$\,\{g_{x_0},R_{x_0},\nabla_{x_0}R,\ldots,\nabla_{x_0}^{r-2}R\}\,$
of covariant derivatives at $\,x_0\,$ of the curvature tensor
$\,R\,$ of $\,g\,$ and vice versa (see \cite{tho}). The main
advantage of using normal tensors is the possibility of expressing
the symmetries of each $\,g_{x_0}^s\,$ without using the other
normal tensors, whereas the symmetries of $\,\nabla_{x_0}^sR\,$
depend on $\,R\,$ (recall the Ricci identities).

\smallskip

\begin{rem}
Using the exact sequence
$$ 0 \longrightarrow N_r  \longrightarrow S^2 T^*_{x_0} X \otimes S^r T^*_{x_0} X  \stackrel{s}{\longrightarrow} T^*_{x_0}X \otimes S^{r+1} T^*_{x_0}X
\longrightarrow 0\,,$$ where $s$ stands for the symmetrization on
the last $(r+1)-$indices, we obtain
$$
\dim N_r\,=\,\binom{n+1}{2}\binom{n+r-1}{r}-n\binom{n+r}{r+1}\,.
$$
\end{rem}

\bigskip

\section{Differential invariants of
metrics}\label{invariantesdiferenciales}

In the remainder of the paper, $\,X\,$ will always be an
$\,n-$dimensional smooth manifold.

Let us denote by $\,J^rM \rightarrow X\,$ the fiber bundle of
$\,r-$jets of semi-Riemannian metrics on $\,X\,$ of fixed
signature $\,(p\,,q)\,$, with $\,n=p+q\,$. Its fiber over a point
$\,x_0\in X\,$ will be denoted $\,J_{x_0}^rM\,$.

Let $\,\mbox{Diff}_{x_0}\,$ be the group of germs of local
diffeomorphisms of $\,X\,$ leaving $\,x_0\,$ fixed, and let
$\,\mbox{Diff}_{x_0}^r\,$ be the Lie group of $\,r-$jets at
$\,x_0\,$ of local diffeomorphisms of $\,X\,$ leaving $\,x_0\,$
fixed. We have the following exact group sequence:
$$0\longrightarrow H_{x_0}^r \longrightarrow \mbox{Diff}_{x_0}
\longrightarrow \mbox{Diff}_{x_0}^r \longrightarrow 0\,,$$
$H_{x_0}^r\,$ being the subgroup of $\,\mbox{Diff}_{x_0}\,$ made
up of those diffeomorphisms whose $\,r-$jet at $\,x_0\,$ coincides
with that of the identity.

The group $\,\mbox{Diff}_{x_0}\,$ acts in an obvious way on
$\,J_{x_0}^rM\,$. Note that the subgroup $\,H_{x_0}^{r+1}\,$ acts
trivially, so the action of $\,\mbox{Diff}_{x_0}\,$ on
$\,J_{x_0}^rM\,$ factors through an action of
$\,\mbox{Diff}_{x_0}^{r+1}\,$.

\smallskip

\begin{defn}
Two $\,r-$jets $\,j_{x_0}^rg\,$,$\,j_{x_0}^r\bar{g}\in
J_{x_0}^rM\,$ are said to be \textbf{equivalent} if there exists a
local diffeomorphism $\,\tau\in \mbox{Diff}_{x_0}\,$ such that
$\,j_{x_0}^r\bar{g}=j_{x_0}^r(\tau^*g)\,$.
\end{defn}

\medskip

Equivalence classes of $\,r-$jets of metrics constitute a ringed
space. To be precise:

\smallskip

\begin{defn}
We call \textbf{moduli space} of $\,r-$jets of semi-Riemannian
metrics of signature $\,(p\,,q)\,$ the quotient ringed space
$$
\mathbb{M}_{p,q}^r:=J_{x_0}^rM/\mbox{Diff}_{x_0} =
J_{x_0}^rM/\mbox{Diff}_{x_0}^{r+1}\,.
$$

In the case of Riemannian metrics, that is $\,p=n\,$, $\,q=0\,$,
the moduli space will be denoted $\,\mathbb{M}_n^r\,$.
\end{defn}

\medskip

It is important to observe that the moduli space depends neither
on the point $\,x_0\,$ nor on the chosen $\,n-$dimensional
manifold:

Given a point $\,\bar{x}_0\,$ in another $\,n-$dimensional
manifold $\,\bar{X}\,$, let us consider an arbitrary
diffeomorphism
$$
\xymatrix{X\supset U_{x_0} \ar[r]^-{\varphi} &
U_{\bar{x}_0}\subset \bar{X}}
$$
between corresponding neighbourhoods of $\,x_0\,$ and
$\,\bar{x}_0\,$, verifying $\,\varphi(x_0)=\bar{x}_0\,$. Such a
diffeomorphism induces an isomorphism of ringed spaces between the
corresponding moduli spaces,
$$
\begin{CD}
J_{\bar{x}_0}^r\bar{M}/\mbox{Diff}_{\bar{x}_0} @= J_{x_0}^rM/\mbox{Diff}_{x_0}\\
[j_{\bar{x}_0}^r\bar{g}] & \longmapsto & [j_{x_0}^r
\varphi^*\bar{g}]\,,
\end{CD}
$$
which is independent of the choice of the diffeomorphism
$\,\varphi\,$. So both moduli spaces are canonically identified.

\smallskip

Let us now consider the quotient morphism
$$
\xymatrix{J_{x_0}^rM \ar[r]^-{\pi} &
J_{x_0}^rM/\mbox{Diff}_{x_0}\,=\,\mathbb{M}_{p,q}^r\,.}
$$

Recall that a function $\,f\,$ defined on an open subset
$\,U\subseteq \mathbb{M}_{p,q}^r\,$ is said to be
\textbf{differentiable} if $\,f\circ\pi\,$ is a smooth function on
$\,\pi^{-1}(U)\,$, that is,
$$
\mathcal{C}^{\infty}(U)
=\mathcal{C}^{\infty}(\pi^{-1}(U))^{\mbox{Diff}_{x_0}}\,.
$$

Every semi-Riemannian metric $\,g\,$ on $\,X\,$ of signature
$\,(p,q)\,$ defines a map
$$
\begin{array}{ccc}
X & \stackrel{m_g}{\longrightarrow} & \mathbb{M}_{p,q}^r\\
x & \longmapsto & [j_x^rg]\,,
\end{array}
$$
which is ``differentiable'', that is, it is a morphism of ringed
spaces.

\smallskip

\begin{defn}\label{definiciondeinvariante}
A \textbf{differential invariant} of order $\,\leq r \,$ of
semi-Riemannian metrics of signature $\,(p,q)\,$ is defined to be
a global differentiable function on $\,\mathbb{M}_{p,q}^r\,$.

Taking into account the ringed space structure of
$\,\mathbb{M}_{p,q}^r\,$, we can simply write:
$$
\{\mbox{Differential invariants of order}\,\,\,\leq
r\}=\mathcal{C}^{\infty}(\mathbb{M}_{p,q}^r)=\mathcal{C}^{\infty}(J_{x_0}^rM)^{\mbox{Diff}_{x_0}}\,.
$$
\end{defn}

\medskip

A differential invariant $\,h:\mathbb{M}_{p,q}^r\rightarrow
\mathbb{R}\,$ associates with every semi-Riemannian metric $\,g\,$
on $\,X\,$ a smooth function on $\,X\,$, denoted by $\,h(g)\,$,
through the formula $\,h(g):=h\circ m_g\,$, that is,
$$h(g)(x)=h([j_x^rg])\,.$$

In any local coordinates, $\,h(g)\,$ is a function smoothly
depending on the coefficients of the metric and their subsequent
partial derivatives up to the order $\,r\,$,
$$
h(g)(x)=h\left( g_{ij}(x), \frac{\partial g_{ij}}{\partial
x_k}(x),\ldots,\frac{\partial^r g_{ij}}{\partial
x_{k_1}\ldots\partial x_{k_r}}(x)\right)\,,
$$
which is equivariant with respect to the action of local
diffeomorphisms,
$$
h(\tau^*g)=\tau^*(h(g))\,.
$$

\medskip

For a discussion on the concept of differential invariant, see
Section \ref{appendix}.

\bigskip

\section{A fundamental lemma}\label{keylemma}

The aim of this section is to prove that there exist a certain
linear finite-dimensional representation $\,V^r\,$ of the
orthogonal group $\,O(p,q)\,$ and an isomorphism of ringed spaces
$$
\mathbb{M}_{p,q}^r\,=\,V^r/\,O(p,q)\,.
$$

This bijection is already known at a set-theoretic level (see
\cite{eps} and also \cite{moduli} for $\,G-$structures which
posses a linear connection). We just add the fact that this
bijection is an isomorphism of ringed spaces.

Let us fix for this entire section a local coordinate system
$\,(z_1,\ldots,z_n)\,$ centred at $\,x_0\,$.

We will denote by $\,\mathcal{N}_{x_0}^r\,$ the smooth submanifold
of $\,J_{x_0}^rM\,$ formed by $\,r-$jets at $\,x_0\,$ of metrics
of signature $\,(p,q)\,$ for which $\,(z_1,\ldots,z_n)\,$ is a
normal coordinate system (that is, Taylor expansions of the
coefficients of such metrics with respect to coordinates
$\,(z_1,\ldots,z_n)\,$ satisfy the equations of the Gauss Lemma up
to the order $\,r\,$).

Consider the subgroup of $\,\mbox{Diff}_{x_0}$
$$
H_{x_0}^1:=\{\tau \in \mbox{Diff}_{x_0}\,:\,
j_{x_0}^1\tau\,=\,j_{x_0}^1(\mbox{Id})\}\,.
$$

Note the following exact group sequence:
$$
0 \longrightarrow H_{x_0}^1 \longrightarrow \mbox{Diff}_{x_0}
\longrightarrow \mbox{Gl}\,(T_{x_0}X) \longrightarrow 0\,,
$$
where the epimorphism $\,\mbox{Diff}_{x_0} \rightarrow
\mbox{Gl}\,(T_{x_0}X)\,$ takes every diffeomorphism to its linear
tangent map at $\,x_0\,$.

\smallskip

\begin{lem}\label{lemaprimero}
There exists an isomorphism of ringed spaces
$$
\xymatrix{\mathcal{N}_{x_0}^r \ar@{=}[r] &
J_{x_0}^rM/H_{x_0}^1\,.}
$$
\end{lem}

\medskip

\begin{proof} Let us start by constructing a smooth section of the
natural inclusion
$$
\mathcal{N}_{x_0}^r \hookrightarrow J_{x_0}^rM\,.
$$

Given a jet metric $\,j_{x_0}^rg\in J_{x_0}^rM\,$, consider a
metric $\,g\,$ representing it. Let
$\,(\bar{z}_1,\ldots,\bar{z}_n)\,$ be the only normal coordinate
system centred at $\,x_0\,$ with respect to $\,g\,$ which
satisfies $\,\mbox{d}_{x_0}\bar{z}_i\,=\,\mbox{d}_{x_0}z_i\,$.

Let $\,\tau\,$ be the local diffeomorphism which transforms one
coordinate system into another: $\,\tau^*(\bar{z}_i)=z_i\,$. The
condition $\,\mbox{d}_{x_0}\bar{z}_i\,=\,\mbox{d}_{x_0}z_i\,$
implies that the linear tangent map of $\,\tau\,$ at $\,x_0\,$ is
the identity, i.e. $\,\tau\in H_{x_0}^1\,$.

As $\,(\bar{z}_1,\ldots,\bar{z}_n)\,$ is a normal coordinate
system for $\,g\,$,
$\,(z_1=\tau^*(\bar{z}_1),\ldots,z_n=\tau^*(\bar{z}_n))\,$ is a
normal coordinate system for $\,\tau^*g\,$; that is,
$\,j_{x_0}^r(\tau^*g)\in \mathcal{N}_{x_0}^r\,$.

Therefore, the section we were looking for is the following map:
$$
\begin{array}{ccc}
J_{x_0}^rM & \stackrel{\varphi}{\longrightarrow} &
\mathcal{N}_{x_0}^r\\
j_{x_0}^rg & \longmapsto & j_{x_0}^r(\tau^*g)\,,
\end{array}
$$
with $\,\tau\,$ depending on $\,g\,$.

Let us now see that $\,\varphi\,$ is constant on each orbit of the
action of $\,H_{x_0}^1\,$. Let $\,j_{x_0}^rg'\,$ be another point
in the same orbit as $\,j_{x_0}^rg\,$, so we can write
$\,g'=\sigma^*g\,$ for some $\,\sigma\in H_{x_0}^1\,$.

Since $\,(\bar{z}_1,\ldots,\bar{z}_n)\,$ are normal coordinates
for $\,g\,$,
$\,(z_1'=\sigma^*(\bar{z}_1),\ldots,z_n'=\sigma^*(\bar{z}_n))\,$
is a normal coordinate system for $\,g'=\sigma^*g\,$. Then
$\,z_i\,=\,\tau^*(\bar{z}_i)\,=\,\tau^*(\sigma^{*^{-1}}(z_i'))\,$,
and, if we apply the definition of $\,\varphi\,$, we get
$$
\varphi(j_{x_0}^rg')\,=\,j_{x_0}^r(\tau^*\sigma^{*^{-1}}g')\,=\,j_{x_0}^r(\tau^*g)\,=\,\varphi(j_{x_0}^rg)\,.
$$

As $\,\varphi\,$ is constant on each orbit of the action of
$\,H_{x_0}^1\,$, it induces, according to the universal quotient
property, a morphism of ringed spaces:
$$
J_{x_0}^rM/H_{x_0}^1 \longrightarrow \mathcal{N}_{x_0}^r\,.
$$

This map is indeed an isomorphism of ringed spaces, because it has
an obvious inverse morphism, which is the following composition:
$$
\mathcal{N}_{x_0}^r \hookrightarrow J_{x_0}^rM \rightarrow
J_{x_0}^rM/H_{x_0}^1\,.
$$
\end{proof}

\smallskip

Let us denote by $\,\mbox{Gl}_n\,$ the general linear group in
dimension $\,n\,$:
$$
\mathrm{Gl}_n:=\,\{\,n\times n\,\,\mathrm{invertible\,\,\,
matrices\,\,\, with\,\,\, coefficients\,\,\,
in}\,\,\mathbb{R}\,\}\,.
$$

Considering every matrix in $\,\mbox{Gl}_n\,$ as a linear
transformation of the coordinate system $\,(z_1,\ldots,z_n)\,$, we
can think of $\,\mbox{Gl}_n\,$ as a subgroup of
$\,\mbox{Diff}_{x_0}\,$.

Via the action of the group $\,\mbox{Diff}_{x_0}\,$ on
$\,J_{x_0}^rM\,$, the subgroup $\,\mbox{Gl}_n\,$, for its part,
acts leaving the submanifold $\,\mathcal{N}_{x_0}^r\,$ stable, and
then we can state the following

\smallskip

\begin{lem}\label{lemasegundo}
There exists an isomorphism of ringed spaces
$$
\xymatrix{\mathcal{N}_{x_0}^r/\,\mathrm{Gl}_n \ar@{=}[r] &
J_{x_0}^rM/\mathrm{Diff}_{x_0}\,=\,\mathbb{M}_{p,q}^r\,.}
$$
\end{lem}

\medskip

\begin{proof} Via the epimorphism
$$
\mathrm{Diff}_{x_0} \longrightarrow
\mathrm{Diff}_{x_0}/\,H_{x_0}^1\,=\,\mathrm{Gl}\,(T_{x_0}X)\,,
$$
the subgroup $\,\mbox{Gl}_n\,$ gets identified with
$\,\mbox{Gl}\,(T_{x_0}X)\,$. Consequently, the subgroups
$\,H_{x_0}^1\,$ and $\,\mbox{Gl}_n\,$ generate
$\,\mbox{Diff}_{x_0}\,$.

If we consider the isomorphism
$$
\xymatrix{\mathcal{N}_{x_0}^r \ar@{=}[r] & J_{x_0}^rM/H_{x_0}^1}
$$
of Lemma \ref{lemaprimero} and take quotient with respect to the
action of $\,\mbox{Gl}_n\,$, we get the desired isomorphism:
$$
\xymatrix{\mathcal{N}_{x_0}^r/\,\mathrm{Gl}_n \ar@{=}[r] &
(J_{x_0}^rM/\,H_{x_0}^1)/\mathrm{Gl}_n\,=\,J_{x_0}^rM/\mathrm{Diff}_{x_0}\,.}
$$
\end{proof}

\smallskip

Let us express the previous result in terms of normal tensors by
using the following

\smallskip

\begin{lem}\label{lematercero}
The map
$$
\xymatrix{\mathcal{N}_{x_0}^r \ar@{=}[r] &
N_0\,\times\,N_2\,\times\,\ldots\,\times\,N_r\,\,\,\,\,,\,\,\,\,\,j_{x_0}^rg\longmapsto
(g_{x_0},g_{x_0}^2,\ldots,g_{x_0}^r)}
$$
is a diffeomorphism.
\end{lem}

\medskip

\begin{proof} The inverse map is defined in the obvious way:

Given $\,(T^0,T^2,\ldots,T^r)\in N_0\times N_2\times \ldots \times
N_r\,$, consider the jet metric $\,j_{x_0}^rg\,$ which in
coordinates $\,(z_1,\ldots,z_n)\,$ is determined by the identities
$$
g_{ij,k_1\ldots k_s}:=\frac{\partial^s g_{ij}}{\partial
z_{k_1}\cdots \partial z_{k_s}}(x_0)\,= \,T_{ijk_1\ldots
k_s}^s\,\,\,,\,\,\, s=0,\ldots,r\,.
$$

The symmetries of tensors $\,T^s\,$ guarantee that the
coefficients $\,g_{ij}\,$ of the metric $\,g\,$ verify the
equations of the Gauss Lemma up to the order $\,r\,$, that is,
$\,j_{x_0}^rg\in \mathcal{N}_{x_0}^r\,$.
\end{proof}

\smallskip

Combining Lemma \ref{lemasegundo} and Lemma \ref{lematercero}, we
obtain an isomorphism of ringed spaces:

$$
\begin{CD}
\mathbb{M}_{p,q}^r\,=\,J_{x_0}^rM /\,\mathrm{Diff}_{x_0} @=
(N_0\times N_2\times\ldots\times
N_r)/\,\mathrm{Gl}\,(T_{x_0}X)\\
[j_{x_0}^rg] & \longmapsto &
[(g_{x_0},g_{x_0}^2,\ldots,g_{x_0}^r)]\,\,.
\end{CD}
$$

Let us now fix a metric $\,g_{x_0}\in N_0\,$ at $\,x_0\,$ and let
us consider the orthogonal group
$\,O(p,q):=O(\,T_{x_0}X,g_{x_0})\,$. As the linear group
$\,\mbox{Gl}(\,T_{x_0}X)\,$ acts transitively on the space of
metrics $\,N_0\,$, and $\,O(p,q)\,$ is the stabilizer subgroup of
$\,g_{x_0}\in N_0\,$, we obtain the following isomorphism:

$$
\xymatrix{(N_0\times N_2\times \ldots \times
N_r)/\,\mathrm{Gl}(\,T_{x_0}X) \ar@{=}[r] & (N_2\times \ldots
\times N_r)/\,O(p,q)\,.}
$$

To sum up, we can state the main result of this section:

\smallskip

\begin{lem} {\bf(Fundamental Lemma)}\label{fundamental}
The moduli space $\,\mathbb{M}_{p,q}^r\,$ is isomorphic to the
quotient space of a linear representation of the orthogonal group
$\,O(p,q)\,$, through the following isomorphism of ringed spaces:
$$
\xymatrix{\mathbb{M}_{p,q}^r \ar@{=}[r] & (N_2\times \ldots \times
N_r) /\,O(p,q)\,.}
$$
\end{lem}

\medskip

This isomorphism takes every class $\,[j_{x_0}^r\bar{g}]\in
\mathbb{M}_{p,q}^r\,$, with $\,\bar{g}_{x_0}=g_{x_0}\,$, to the
sequence of normal tensors
$\,[(\bar{g}_{x_0}^2,\ldots,\bar{g}_{x_0}^r)]\in (N_2\times \ldots
\times N_r)/\,O(p,q)\,.$

\bigskip

\section{Structure of the moduli spaces}\label{structure}

Let $\,V\,$ be a finite-dimensional linear representation of a
reductive Lie group $\,G\,$. The $\,\mathbb{R}-$algebra of
$\,G-$invariant polynomials on $\,V\,$ is finitely generated
(Hilbert-Nagata theorem, see \cite{fog}). Let $\,p_1,\ldots,p_s\,$
be a finite set of generators for that algebra; by a result of
Luna \cite{luna}, every smooth $\,G-$invariant function $\,f\,$ on
$\,V\,$ can be written as $\,f=F(p_1,\ldots,p_s)\,$, for some
smooth function $\,F\in \mathcal{C}^{\infty}(\mathbb{R}^s)\,$.

\smallskip

\begin{thm} {\bf(Finiteness of differential invariants,
\cite{masval})}\label{finitud} There exists a finite number
$\,p_1,\ldots,p_s\in \mathcal{C}^{\infty}(\mathbb{M}_{p,q}^r)\,$
of differential invariants of order $\,\leq r\,$ such that any
other differential invariant $\,f\,$ of order $\,\leq r\,$ is a
smooth function of the former ones, i.e.
$\,f=F(p_1,\ldots,p_s)\,$, for a certain $\,F\in
\mathcal{C}^{\infty}(\mathbb{R}^s)\,$.
\end{thm}

\medskip

\begin{proof} By the Fundamental Lemma (\ref{fundamental}),
$$
\mathcal{C}^{\infty}(\mathbb{M}_{p,q}^r)\,=\,\mathcal{C}^{\infty}(N_2\times
\ldots \times N_r)^{O(p,q)}\,,
$$
and we can conclude by applying the above theorem by Luna to the
linear representation $\,N_2\times\ldots\times N_r\,$ of the
orthogonal group $\,O(p,q)\,$.
\end{proof}

\smallskip

\begin{rem} Using the theory of invariants for the orthogonal
group and the fact that the sequence of normal tensors
$\,\{g_{x_0},g_{x_0}^2,g_{x_0}^3,\ldots,g_{x_0}^r\}\,$ is
equivalent to the sequence
$\,\{g_{x_0},R_{x_0},\nabla_{x_0}R,\ldots,\nabla_{x_0}^{r-2}R\}\,$,
it can be proved that the generators $\,p_1,\ldots,p_s\,$ of
Theorem \ref{finitud} can be chosen to be \textbf{Weyl
invariants}, that is, scalar quantities constructed from the
sequence
$\,\{g_{x_0},R_{x_0},\nabla_{x_0}R,\ldots,\nabla_{x_0}^{r-2}R\}\,$
by reiteration of the following operations: tensor products,
raising and lowering indices, and contractions.
\end{rem}

\medskip

\begin{thm}\label{clasifica} In the Riemannian case, differential invariants of
order $\,\leq r\,$ separate points in the moduli space
$\,\mathbb{M}_n^r\,$.

Consequently, differential invariants of order $\,\leq r\,$
classify $\,r-$jets of Riemannian metrics (at a point).
\end{thm}

\medskip

\begin{proof} For positive definite metrics, the orthogonal group
$\,O(n)\,$ is compact. It is a well-known fact that, if $\,V\,$ is
a linear representation of a compact Lie group $\,G\,$, then
smooth $\,G-$invariant functions on $\,V\,$ separate the orbits of
the action of $\,G\,$, or, in other words, the algebra
$\,\mathcal{C}^{\infty}(V/G)\,$ separates the points in $\,V/G\,$.

Using this, together with the Fundamental Lemma, we conclude our
proof.
\end{proof}

\smallskip

Neither assertion in Theorem \ref{clasifica} is valid for
semi-Riemannian metrics. See Note in Subsection \ref{ejemplos}\,
for a counterexample. For such metrics, moduli spaces
$\,\mathbb{M}_{p,q}^r\,$ are generally pathological in a
topological sense, since they have non-closed points (they are not
$\,T_1\,$ topological spaces).

\smallskip

In the Riemannian case, Schwarz Theorem \ref{schwarz} and the
Fundamental Lemma directly provide the following

\smallskip

\begin{thm} In the Riemannian case, moduli spaces
$\,\mathbb{M}_n^r\,$ are differentiable spaces.
\end{thm}

\medskip

More precisely: Let $\,p_1,\ldots,p_s\,$ be the basis of
differential invariants of order $\,\leq r\,$ mentioned in Theorem
\ref{finitud}. These invariants induce an isomorphism of
differentiable spaces
$$
\xymatrix{(p_1,\ldots,p_s)\,:\,\mathbb{M}_n^r \ar@{=}[r] &
Z\,\subseteq\,\mathbb{R}^s\,,}
$$
$Z\,$ being a closed subspace of $\,\mathbb{R}^s\,$.

Although the differentiable space $\,\mathbb{M}_n^r\,$ is not in
general a smooth manifold, its structure is not so deficient as it
could seem at first sight, since we are going to prove that it
admits a finite stratification by certain smooth submanifolds.

\smallskip

\begin{defn} Let us consider $\,V_n=\mathbb{R}^n\,$ endowed with
its standard inner product $\,\delta\,$, and the corresponding
orthogonal group $\,O(n):=O(V_n,\delta)\,$. We will denote by
$\,\mathcal{T}\,$ the set of conjugacy classes of closed subgroups
in $\,O(n)\,$.

Given another $\,n-$dimensional vector space $\,\bar{V}_n\,$ with
an inner product $\,\bar{\delta}\,$, we can also consider the set
$\,\bar{\mathcal{T}}\,$ of conjugacy classes of closed subgroups
in $\,O(\bar{V}_n,\bar{\delta})\,$.

Observe that there exists a canonical identification
$$
\mathcal{T} \longrightarrow \bar{\mathcal{T}}\,\,\,,\,\,\, [H]
\longmapsto [\varphi \circ H \circ \varphi^{-1}]\,,
$$
where $\,\varphi\,$ stands for any isometry
$\,\varphi:V_n\rightarrow \bar{V}_n\,$.

As the identification is canonical (i.e. it does not depend on the
choice of the isometry $\,\varphi\,$), from now on we will suppose
that the set $\,\mathcal{T}\,$ is just ``the same'' for every pair
$\,(\bar{V}_n,\bar{\delta})\,$.

Note that $\,\mathcal{T}\,$ possesses a partial order relation:
$\,[H] \leq [H']\,$, if there exist some representatives $\,H\,$
and $\,H'\,$ of $\,[H]\,$ and $\,[H']\,$ respectively, such that
$\,H\subseteq H'\,$.
\end{defn}

\medskip

\begin{defn} The \textbf{group of automorphisms} of a Riemannian
jet metric $\,j_{x_0}^rg\,$ is defined to be the stabilizer
subgroup $\,\mbox{Aut}(j_{x_0}^rg)\subseteq
\mbox{Diff}_{x_0}^{r+1}\,$ of $\,j_{x_0}^rg\,$:
$$
\mathrm{Aut}(j_{x_0}^rg):=\{j_{x_0}^{r+1}\tau\in
\mathrm{Diff}_{x_0}^{r+1}\,:\,j_{x_0}^r(\tau^*g)=j_{x_0}^rg\}\,.
$$
\end{defn}

\medskip

Given $\,\tau\in \mbox{Diff}_{x_0}\,$, let us denote by
$\,\tau_{*,x_0}:T_{x_0}X\rightarrow T_{x_0}X\,$ the linear tangent
map of $\,\tau\,$ at $\,x_0\,$.

\smallskip

\begin{lem} The group morphism
$$
\begin{array}{ccc}
\mathrm{Aut}(j_{x_0}^rg) & \longrightarrow &
O(T_{x_0}X,g_{x_0})\simeq O(n)\\
j_{x_0}^{r+1}\tau & \longmapsto & \tau_{*,x_0}
\end{array}
$$
is injective.
\end{lem}

\medskip

\begin{proof} For any $\,\tau\in \mbox{Diff}_{x_0}\,$ and any
metric $\,g\,$ on $\,X\,$ we have the following commutative
diagram of local diffeomorphisms:
$$
\xymatrix{T_{x_0}X \ar[r]^-{\mathrm{exp}_{\tau^*\!g}}
\ar[d]_-{\tau_*} & X \ar[d]^-{\tau} \\
T_{x_0}X \ar[r]^-{\mathrm{exp}_g} & X}
$$

If $j_{x_0}^{r+1}\tau\in \mbox{Aut}(j_{x_0}^rg)\,$, that is,
$\,j_{x_0}^r(\tau^*g)=j_{x_0}^rg\,$, then
$\,j_0^{r+1}(\mbox{exp}_{\tau^*g})=j_0^{r+1}(\mbox{exp}_g)\,$
because of Proposition \ref{exponencial}.

Now, taking $\,(r+1)-$jets in the above diagram, we obtain:
$$
j_{x_0}^{r+1}\tau\,=\,j_0^{r+1}(\mathrm{exp}_g)\,\circ\,j_0^{r+1}\tau_*\,\circ\,j_{x_0}^{r+1}(\mathrm{exp}_g^{-1})\,,
$$
hence $\,j_{x_0}^{r+1}\tau\,$ is determined by its linear part
$\,\tau_*\,$.
\end{proof}

\smallskip

By the previous lemma, the group $\,\mbox{Aut}(j_{x_0}^rg)\,$ can
be viewed as a subgroup (determined up to conjugacy) of the
orthogonal group $\,O(n)\,$.

\smallskip

\begin{defn} The \textbf{type map} is defined to be the map
$$
t:\mathbb{M}_n^r \longrightarrow \mathcal{T} \,\,\,,\,\,\,
[j_{x_0}^rg] \longmapsto [\mathrm{Aut}(j_{x_0}^rg)]\,.
$$

For each $\,[H]\in \mathcal{T}\,$, the \textbf{stratum of type}
$\,[H]\,$ is said to be the subset $\,S_{[H]}\subseteq
\mathbb{M}_n^r\,$ of those points of type $\,[H]\,$.
\end{defn}

\medskip

\begin{thm}\label{estratifica} {\bf(Stratification of the moduli space)}
The type map $\,t:\mathbb{M}_n^r\rightarrow \mathcal{T}\,$
verifies the following properties:

\begin{enumerate}
    \item $\,t\,$ takes a finite number of values $\,[H_0],\ldots,[H_k]\,$,
    one of which, say $\,[H_0]\,$, is minimum.

    \item Semicontinuity: For every type $\,[H]\in \mathcal{T}\,$, the
    set of points in $\,\mathbb{M}_n^r\,$ of type $\,\leq [H]\,$ is an
    open subset of $\,\mathbb{M}_n^r\,$. In particular, every
    stratum $\,S_{[H_i]}\,$ is a locally closed subspace of
    $\,\mathbb{M}_n^r\,$.

    \item Every stratum $\,S_{[H_i]}\,$ is a smooth submanifold of
    $\,\mathbb{M}_n^r\,$.

    \item The (also called generic) stratum $\,S_{[H_0]}\,$ of minimum
    type is a dense connected open subset of $\,\mathbb{M}_n^r\,$.
  \end{enumerate}

\end{thm}

\medskip

\begin{proof} Fix a positive definite metric $\,g_{x_0}\,$ on
$\,T_{x_0}X\,$ and denote by $\,O(n)\,$ its orthogonal group. The
Fundamental Lemma \ref{fundamental} tells us that there exists an
isomorphism
$$
\xymatrix{\mathbb{M}_n^r \ar@{=}[r] & (N_2\times \ldots \times
N_r)/\,O(n)\,.}
$$

This isomorphism takes every class
$\,[j_{x_0}^r\bar{g}]\in\mathbb{M}_n^r\,$, with
$\,\bar{g}_{x_0}=g_{x_0}\,$, to the sequence of normal tensors
$\,[\bar{g}_{x_0}^2,\ldots,\bar{g}_{x_0}^r]\in
(N_2\times\ldots\times N_r)/\,O(n)\,$.

Let us check that the subgroup $\,\mbox{Aut}(j_{x_0}^r\bar{g})
\hookrightarrow O(n)\,\,,\,\,j_{x_0}^{r+1}\tau \mapsto \tau_*\,$,
coincides with the subgroup
$$
\mathrm{Aut}(\bar{g}_{x_0}^2,\ldots,\bar{g}_{x_0}^r):=\{\sigma\in
O(n):
\sigma^*(\bar{g}_{x_0}^k)\,=\,\bar{g}_{x_0}^k\,,\,\forall\,k\leq
r\}\,.
$$

It is clear that if an automorphism $\,j_{x_0}^{r+1}\tau\,$ leaves
$\,j_{x_0}^r\bar{g}\,$ fixed, then the sequence of its normal
tensors must also remain fixed by the automorphism:
$\,\tau^*(\bar{g}_{x_0}^k)\,=\,\bar{g}_{x_0}^k\,$.

\smallskip

Reciprocally, given an automorphism $\,\sigma:T_{x_0}X\rightarrow
T_{x_0}X\,$ of the sequence of normal tensors
$\,(\bar{g}_{x_0}^2,\ldots,\bar{g}_{x_0}^r)\,$, let us consider a
normal coordinate system $\,z_1,\ldots,z_n\,$ for $\,\bar{g}\,$ at
$\,x_0\,$.

Via the identification provided by the exponential map
$\,\mbox{exp}_g:T_{x_0}X\rightarrow X\,$, the map $\,\sigma\,$ can
be viewed as a diffeomorphism of $\,X\,$ (a linear transformation
of normal coordinates).

In normal coordinates, the expression of the normal tensor
$\,\bar{g}_{x_0}^k\,$ corresponds to the expression of the
homogeneous part of degree $\,k\,$ of the jet metric
$\,j_{x_0}^r\bar{g}\,$. Hence it is an immediate consequence that
the linear transformation $\,\sigma\,$ leaves
$\,j_{x_0}^r\bar{g}\,$ fixed, i.e.
$\,j_{x_0}^{r+1}\sigma\in\mbox{Aut}(j_{x_0}^r\bar{g})\,$.

The identity
$\,\mbox{Aut}(j_{x_0}^r\bar{g})\,=\,\mbox{Aut}(\bar{g}_{x_0}^2,\ldots,\bar{g}_{x_0}^r)\,$
implies that the following diagram is commutative:

$$
\begin{CD}
\mathbb{M}_n^r @>t>> \mathcal{T}\\
@| @|\\
(N_2\times\ldots\times N_r)/\,O(n) @>t>> \mathcal{T}\\
[\bar{g}_{x_0}^2,\ldots,\bar{g}_{x_0}^r] & \longmapsto &
[\mathrm{Aut}(\bar{g}_{x_0}^2,\ldots,\bar{g}_{x_0}^r)]\,.
\end{CD}
$$

Therefore, our theorem has come down to the case of a linear
representation $\,V(=N_2\times\ldots\times N_r)\,$ of a compact
Lie group $\,G(=O(n))\,$ and the corresponding type map:
$$
\begin{array}{ccl}
V/G & \stackrel{t}{\longrightarrow} &
\mathcal{T}=\{\mathrm{conjugacy\,\,classes\,\,of\,\,closed\,\,subgroups\,\,of\,\,}G\}\\

[v] & \longmapsto &
[\mathrm{Stabilizer\,\,subgroup\,\,of\,\,}v]\,.
\end{array}
$$

For this type map, the analogous properties to 1 -- 4 in the
statement are well known (see \cite{bourbaki}, Chap. IX, \S\,9,
Th.\! 2 and Exer.\! 9).
\end{proof}

\smallskip

\begin{rem} Except for trivial cases, the generic stratum has type
$\,H_0=\{0\}\,$.
\end{rem}

\medskip

\begin{rem} The dimension of the moduli space $\,\mathbb{M}_n^r\,$
(or rather that of its generic stratum) can be deduced directly
from the Fundamental Lemma and the formulae giving the dimensions
of spaces $\,N_r\,$ of normal tensors which were presented in
Section \ref{preliminares}.

The result (due, in a different language, to J. Mu{\~n}oz and A.
Vald\'{e}s, \cite{masval2}) is as follows:
\[
\dim \mathbb{M}_n^0\,=\,\dim
\mathbb{M}_n^1\,=\,0\,\,\,,\,\,\,\forall\, n\geq 1\,;
\]
\[
\dim \mathbb{M}_1^r\,=\,0\,\,\,,\,\,\,\forall\, r\geq 0\,;
\]
\[
\dim \mathbb{M}_2^2\,=\,1\,\,\,\,\,\,,\,\,\,\,\,\, \dim
\mathbb{M}_2^r=\frac{1}{2}(r+1)(r-2)\,\,\,,\,\,\,\forall\, r\geq
3\,;
\]
\[
\dim
\mathbb{M}_n^r\,=\,n+\frac{(r-1)n^2-(r+1)n}{2(r+1)}\binom{n+r}{r}\,\,\,,\,\,\,\forall\,
n\geq 3\,,\,r\geq 2\,.
\]
\end{rem}

\bigskip

\section{Moduli spaces in dimension $\,n=2\,$}\label{superficies}

\subsection{Stratification}

We are going to determine the stratification of moduli spaces
$\,\mathbb{M}_2^r\,$ of $\,r-$jets of Riemannian metrics in
dimension $\,n=2\,$.

Let us consider the vector space $\,\mathbb{R}^2=\mathbb{C}\,$,
endowed with the standard Euclidean metric, and its corresponding
orthogonal group $\,O(2)\,$. We will denote by $\,(x,y)\,$ the
Cartesian coordinates and by $\,z=x+iy\,$ the complex coordinate.

Let us denote by $\,\sigma_m:\mathbb{C}\rightarrow \mathbb{C}\,$
the rotation of angle $\,2\pi/m\,$ (that is,
$\,\sigma_m(z)=\varepsilon_m z\,$, with
$\,\varepsilon_m=\cos(2\pi/m)\,+\,i\sin(2\pi/m)\,$ a primitive
$\,m$th root of unity) and by $\,\tau:\mathbb{C}\rightarrow
\mathbb{C}\,,\,\tau(z)=\bar{z}\,$ the complex conjugation.

The only (up to conjucagy) closed subgroups of $\,O(2)\,$ are the
following ones:
\[
SO(2):=\{\varphi\in O(2):\det
\varphi=1\}\,\,\,\,\mathrm{(special\,\, orthogonal\,\,
group)}\,,\]
\[
K_m:=<\sigma_m>\,\,\,\,\,\,\,\,\,\,\,\,\mathrm{(group\,\, of\,\,
rotations\,\, of\,\, order}\,\,
m\,\mathrm{)}\,\,\,\,\,\,\,\mathrm{(}m\geq 1\mathrm{)}\,,
\]
\[
D_m:=<\sigma_m,\tau>\,\,\,\,\,\,\,\,\,\,\,\,\mathrm{(dihedral\,\,
group\,\, of\,\, order}\,\,
2m\,\mathrm{)}\,\,\,\,\,\,\,\mathrm{(}m\geq 1\mathrm{)}\,,
\]
and $\,O(2)\,$ itself. All these subgroups are normal but the
dihedral $\,D_m\,$.

The subgroup $\,SO(2)\,$ of rotations is identified with the
multiplicative group $\,S_1\subset \mathbb{C}\,$ of complex
numbers of modulus 1,
$$
\begin{CD}
S_1 @= SO(2)\\
\alpha & \longmapsto & \rho_{\alpha}\,\,\,\,\,\,\,, & &
\,\,\,\,\,\,\,\,\rho_{\alpha}(z):=\alpha z\,.
\end{CD}
$$
Besides, every element in $\,O(2)\,$ is either $\,\rho_{\alpha}\,$
or $\,\tau\rho_{\alpha}\,$, for some $\alpha\in S_1\,$.

The action of $\,O(2)\,$ on $\,\mathbb{R}^2\,$ induces an action
on the algebra $\,\mathbb{R}[x,y]\,$ of the polynomials on
$\,\mathbb{R}^2\,$, to be more specific: $\,\varphi \cdot
P(x,y):=P(\varphi^{-1}(x,y))\,$.

The following lemma provides us with the list of all invariant
polynomials with respect to each of the subgroups of $\,O(2)\,$
above mentioned:

\smallskip

\begin{lem}\label{polinvariantes} The following identities hold:
\begin{enumerate}

  \item
  $\,\mathbb{R}[x,y]^{K_m}\,=\,\mathbb{R}[x^2+y^2,p_m(x,y),q_m(x,y)]\,\,,$

  \item
  $\,\mathbb{R}[x,y]^{D_m}\,=\,\mathbb{R}[x^2+y^2,p_m(x,y)]\,\,,$

  \item
  $\,\mathbb{R}[x,y]^{O(2)}\,=\,\mathbb{R}[x,y]^{SO(2)}\,=\,\mathbb{R}[x^2+y^2]\,\,,$
\end{enumerate}
with $\,p_m(x,y)=\mathrm{Re}(z^m)\,$ and
$\,q_m(x,y)=\mathrm{Im}(z^m)\,$.
\end{lem}

\medskip

\begin{proof}\,\,\,1.\,\, Let us consider the algebra
of polynomials on $\,\mathbb{R}^2\,$ with complex coefficients,
$$
\mathbb{C}[x,y]\,=\,\mathbb{C}[z,\bar{z}]\,=\,\bigoplus_{ab}\mathbb{C}z^{a}\bar{z}^b\,.
$$
Every summand is stable under the action of $\,K_m\,$, since
$$
\sigma_m\cdot(z^{a}\bar{z}^b)\,=\,\frac{1}{\varepsilon_m^{a}\bar{\varepsilon}_m^b}z^{a}\bar{z}^b\,=\,\varepsilon_m^{b-a}z^{a}\bar{z}^b\,.
$$

This formula also tells us that the monomial $\,z^{a}\bar{z}^b\,$
is invariant by $\,K_m\,$ if and only if
$b-a\,\equiv\,0\,\,\mbox{mod}\,m\,$, that is, $\,b-a\,=\,\pm km\,$
for some $\,k\in \mathbb{N}\,$. Then invariant monomials are of
the form
$$
z^{a}\bar{z}^b\,=\,(z\bar{z})^{a}\bar{z}^{km}\,\,\,\,\,\mathrm{or}\,\,\,\,\,z^{a}\bar{z}^b\,=\,(z\bar{z})^b
z^{km}\,\,,
$$
whence
$$
\mathbb{C}[x,y]^{K_m}\,=\,\mathbb{C}[z\bar{z},z^m,\bar{z}^m]\,.
$$

As $\,z\bar{z}\,=\,x^2+y^2\,$, $\,z^m+\bar{z}^m\,=\,2p_m(x,y)\,$
and $\,z^m-\bar{z}^m\,=\,2iq_m(x,y)\,$, we can conclude that
$$
\mathbb{C}[x,y]^{K_m}\,=\,\mathbb{C}[x^2+y^2,p_m(x,y),q_m(x,y)]\,,
$$
and particularly,
$$
\mathbb{R}[x,y]^{K_m}\,=\,\mathbb{R}[x^2+y^2,p_m(x,y),q_m(x,y)]\,.
$$

\smallskip

2.\,\, As $\,D_m=<K_m,\tau>\,$, we get
$$
\mathbb{C}[x,y]^{D_m}\,=\,(\mathbb{C}[x,y]^{K_m})^{<\tau>}\,=\,\mathbb{C}[z\bar{z},z^m,\bar{z}^m]^{<\tau>}
$$
$$
=\,\left[ \left( \bigoplus_k \mathbb{C}[z\bar{z}]z^{km}\right)
\oplus \left( \bigoplus_k \mathbb{C}[z\bar{z}]\bar{z}^{km}\right)
\right]^{<\tau>}
$$
(as $\,\tau\cdot z=\bar{z}\,$ and $\,\tau\cdot \bar{z}=z\,$)
$$
=\,\bigoplus_k
\mathbb{C}[z\bar{z}](z^{km}+\bar{z}^{km})\,=\,\mathbb{C}[z\bar{z},z^m+\bar{z}^m]\,=\,\mathbb{C}[x^2+y^2,p_m(x,y)]\,,
$$
and, in particular,
$$
\mathbb{R}[x,y]^{D_m}\,=\,\mathbb{R}[x^2+y^2,p_m(x,y)]\,.
$$

\smallskip

3.\,\, Every summand in the decomposition
$$
\mathbb{C}[z,\bar{z}]\,=\,\bigoplus_{ab}\mathbb{C}z^{a}\bar{z}^b
$$
is stable under the action of $\,SO(2)\,$, since for every
$\,\rho_{\alpha}\in SO(2)\,$ it is satisfied:
$$
\rho_{\alpha}\cdot(z^{a}\bar{z}^b)\,=\,\frac{1}{\alpha^{a}\bar{\alpha}^b}z^{a}\bar{z}^b\,.
$$

Moreover, this formula assures us that the only monomials
$\,z^{a}\bar{z}^b\,$ which are $\,SO(2)-$in\-var\-iant are those
verifying $\,a=b\,$. Then,
$$
\mathbb{C}[x,y]^{SO(2)}\,=\,\mathbb{C}[z,\bar{z}]^{SO(2)}\,=\,\mathbb{C}[z\bar{z}]\,=\,\mathbb{C}[x^2+y^2]\,,
$$
whence
$$
\mathbb{R}[x,y]^{SO(2)}\,=\,\mathbb{R}[x^2+y^2]\,.
$$

Finally, this identity tells us that $\,SO(2)-$invariant
polynomials are $\,O(2)-$invariant too, so the obvious inclusion
$\,\mathbb{R}[x,y]^{O(2)}\,\subseteq\,\mathbb{R}[x,y]^{SO(2)}\,$
is indeed an equality.
\end{proof}

\smallskip

\begin{cor}\label{estabilizadores} With the same notations used in the previous lemma, it
is verified:

1. $D_m\,$ is the stabilizer subgroup of the polynomial
$\,p_m(x,y)\,$, and there exists no polynomial in
$\,\mathbb{R}[x,y]\,$ of degree $\,<m\,$ whose stabilizer subgroup
is $\,D_m\,$.

2. $K_m$ \!($m\geq 2$) is the stabilizer subgroup of the
polynomial $p_m(x,y)+(x^2+y^2)q_m(x,y)\,$, and there exists no
polynomial in $\,\mathbb{R}[x,y]\,$ of degree $\,<m+2\,$ whose
stabilizer subgroup is $\,K_m\,$.

3. $K_1=\{\mathrm{Id}\}\,$ is the stabilizer subgroup of the
polynomial $\,x+xy\,$, and there exists no polynomial in
$\,\mathbb{R}[x,y]\,$ of degree $\,<2\,$ whose stabilizer subgroup
is $\,K_1\,$.
\end{cor}

\medskip

\begin{proof}\,\,\,1.\,\, Using that every element in $\,O(2)\,$
is either of the form $\,\rho_{\alpha}\,$ or of the form
$\,\rho_{\alpha}\circ \tau\,$, it is a matter of routine to check
that the stabilizer subgroup of the polynomial
$\,p_m(x,y)=\mbox{Re}(z^m)\,$ is $\,D_m\,$.

If there were another polynomial $\,\bar{p}(x,y)\,$ of degree
$\,<m\,$ with the same property, $\,\bar{p}(x,y)\,$ should be a
power of $\,x^2+y^2\,$, because of Lemma \ref{polinvariantes}(2),
and in that case its stabilizer subgroup would be the whole
$\,O(2)\,$, against our hipothesis.

\smallskip

2.\,\, According to Lemma \ref{polinvariantes} (1), every
$\,K_m-$invariant polynomial of degree $\,\leq m\,$ is of the form
$\,\lambda p_m(x,y)+\mu q_m(x,y)\,$ (up to addition of a power of
$\,x^2+y^2\,$). However, a polynomial of such a form does not have
$\,K_m\,$ as its stabilizer subgroup, but a larger dihedral group:
after multiplying by a scalar, we can indeed assume
$\,\lambda^2+\mu^2=1\,$; if $\,\alpha=\lambda-i\mu\,$, then
$$
\lambda p_m(x,y)+\mu q_m(x,y)\,=\,\mathrm{Re}(\alpha
z^m)\,=\,\mathrm{Re}((\beta z)^m)
$$
(with $\,\beta^m=\alpha\,$)
$$
=\rho_{\beta^{-1}}\cdot
\mathrm{Re}(z^m)\,=\,\rho_{\beta^{-1}}\cdot p_m(x,y)\,,
$$
whose stabilizer subgroup is the dihedral group
$\,\rho_{\beta^{-1}}\cdot D_m \cdot \rho_{\beta}\,$, which is
conjugate to the stabilizer subgroup $\,D_m\,$ of $\,p_m(x,y)\,$.
(In particular, taking $\,\lambda=0\,$, $\,\mu=-1\,$, we get that
the stabilizer subgroup of $\,q_m(x,y)\,$ is
$\,\rho_{\beta^{-1}}\cdot D_m \cdot \rho_{\beta}\,$, for
$\,\beta^m=i\,$).

As no polynomial of degree $\,\leq m\,$ has the desired stabilizer
subgroup $\,K_m\,$, and there are not any $\,K_m-$invariant
polynomials of degree $\,m+1\,$ (up to a power of $\,x^2+y^2\,$),
the following degree to be considered is $\,m+2\,$. The stabilizer
subgroup of the polynomial $\,p_m(x,y)+(x^2+y^2)q_m(x,y)\,$, of
degree $\,m+2\,$, is the intersection of the stabilizer subgroups
of its two homogeneous components, $\,p_m(x,y)\,$ and
$\,(x^2+y^2)q_m(x,y)\,$, that is,
$$
D_m \cap (\rho_{\beta^{-1}}\cdot D_m \cdot
\rho_{\beta})\,=\,K_m\,\,\,\,\,\,\,(\beta^m=i)\,.
$$

\smallskip

3.\,\, This case is trivial.
\end{proof}

\smallskip

\begin{thm} The strata in the moduli space $\,\mathbb{M}_2^r\,$
correspond exactly to the following types:
$\,[O(2)]\,,\,[D_1]\,,\ldots,\,[D_{r-2}]\,,\,[K_1]\,,\,\ldots\,,\,[K_{r-4}]\,$.
(And also $\,[K_1]\,$, if $\,r = 4\,$.)
\end{thm}

\medskip

\begin{proof} It is a classical result (see \cite{eps}) that in
dimension 2\, every Riemannian metric can be written in normal
coordinates $\,(x,y)\,$ (in a unique way up to an orthogonal
transformation) as follows:
$$
g\,=\,\mathrm{d}x^2+\mathrm{d}y^2+h(x,y)(y\mathrm{d}x-x\mathrm{d}y)^2\,,
$$
for some smooth function $\,h(x,y)\,$.

Observe that the stabilizer subgroup of $\,O(2)\,$ for the jet
$\,j_0^kh\,$ is the same as that for $\,j_0^{k+2}g\,$.

If we take $\,h(x,y)=0\,$, we get a metric (the Euclidean one,
i.e. $\,g=\mbox{d}x^2+\mbox{d}y^2\,$) whose group of automorphisms
(for any jet order) is $\,O(2)\,$.

Choosing $\,h(x,y)=p_m(x,y)\,$, we obtain an $\,r-$jet metric
(with $\,r\geq m+2\,$) whose stabilizer subgroup is $\,D_m\,$,
because of Corollary \ref{estabilizadores} (1).

If we choose $\,h(x,y)=p_m(x,y)+(x^2+y^2)q_m(x,y)\,$, we get an
$\,r-$jet metric (with $\,r\geq m+4\,$) whose stabilizer subgroup
is $\,K_m\,$, by Corollary \ref{estabilizadores} (2).

If we make $\,h(x,y)=x+xy\,$, then we get an $\,r-$jet metric
(with $\,r\geq 4\,$) whose stabilizer subgroup is $\,K_1\,$,
according to Corollary \ref{estabilizadores} (3).

Finally, let us note that no $\,r-$jet metric can have $\,SO(2)\,$
as its stabilizer subgroup, since such a metric would correspond
to a jet function $\,j_0^{r-2}h\,$ whose stabilizer subgroup
should be $\,SO(2)\,$, which is impossible, because, by Lemma
\ref{polinvariantes} (3), every $\,SO(2)-$invariant polynomial is
also $\,O(2)-$invariant.
\end{proof}

\smallskip

\begin{cor} Every closed subgroup of $\,O(2)\,$, except for
$\,SO(2)\,$, is the group of automorphisms of a jet metric
$\,j_0^rg\,$ on $\,\mathbb{R}^2\,$ for some order $\,r\,$.
\end{cor}

\medskip

\begin{cor} The number of strata in $\,\mathbb{M}_2^r\,$ is:
$$
\mathrm{Number\,\,of\,\,strata\,\,in\,\,}\mathbb{M}_2^r\,=\left \{
\begin{array}{lcl}
1 & & \mathrm{for\,\,} r=0,1,2\\
2 & & \mathrm{for\,\,} r=3\\
4 & & \mathrm{for\,\,} r=4\\
2r-5 & & \mathrm{for\,\,} r\geq 5
\end{array}
\right .
$$
\end{cor}

\bigskip

\subsection{Examples}\label{ejemplos}

Now we describe, without proofs, low order jets in dimension
$\,n=2\,$.

For order $\,r=0,1\,$ (and in any dimension $\,n\,$) moduli spaces
$\,\mathbb{M}_n^r\,$ come down to a single point.

\smallskip

\begin{center}
\textbf{Case $\,r=2\,$\,.}
\end{center}

The moduli space is a line:
$$
\xymatrix{\mathbb{M}_2^2 \ar@{=}[r] & \mathbb{R}\,\,\,\,,\,\,\,\,
[j_{x_0}^2g] \ar@{|->}[r] & K_g(x_0)\,\,.}
$$
In other words, the curvature classifies $\,2-$jets of Riemannian
metrics in dimension $\,n=2\,$.

In this case there is just one stratum, the generic one, whose
type is $\,[O(2)]\,$.

\smallskip

\begin{center}
\textbf{Case $\,r=3\,$\,.}
\end{center}

The moduli space is a closed semiplane:
$$
\xymatrix{\mathbb{M}_2^3 \ar@{=}[r] & \mathbb{R}\times
[0,+\infty)\,\,\,\,,\,\,\,\,[j_{x_0}^3g] \ar@{|->}[r] &
(K_g(x_0)\,,|\mathrm{grad}_{x_0}K_g|^2)\,\,.}
$$
That is to say, the curvature and the square of the modulus of the
gradient of the curvature classify $\,3-$jet metrics in dimension
$\,n=2\,$.

Now we have two different strata:

The generic stratum $\,S_{[D_1]}\,=\,\mathbb{R}\times
(0,+\infty)\,$, with type $\,[D_1]\,$. This stratum is the set of
all classes of jets $\,j_{x_0}^3g\,$ verifying
$\,\mbox{grad}_{x_0}K_g \neq 0\,$ (in this case, the group of
automorphisms is the group of order 2\, generated by the
reflection across the vector $\,\mbox{grad}_{x_0}K_g\,$).

The non-generic stratum $\,S_{[O(2)]}\,=\,\mathbb{R}\times
\{0\}\,$, with type $\,[O(2)]\,$, is the set of all classes of
jets $\,j_{x_0}^3g\,$ verifying $\,\mbox{grad}_{x_0}K_g\,=\,0\,$
(which are invariant with respect to every orthogonal
transformation of normal coordinates).

\smallskip

{\bf Note:} If we consider metrics of signature $\,(+,-)\,$,
instead of Riemannian metrics, then the map
$$
\xymatrix{\mathbb{M}_2^3 \ar[r] & \mathbb{R}\times
[0,+\infty)\,\,\,\,,\,\,\,\,[j_{x_0}^3g] \ar@{|->}[r] &
(K_g(x_0)\,,|\mathrm{grad}_{x_0}K_g|^2)\,\,.}
$$
is not injective, that is, differential invariants do not classify
$\,3-$jet metrics of signature $\,(+,-)\,$. To illustrate this,
consider two metrics $\,g\,$,$\,\bar{g}\,$ of signature
$\,(+,-)\,$, such that $\,K_g(x_0)\,=\,K_{\bar{g}}(x_0)\,$,
$\,\mbox{grad}_{x_0}K_g\,=\,0\,$ and
$\,\mbox{grad}_{x_0}K_{\bar{g}}\,$ is a non-zero isotropic vector
with respect to $\,\bar{g}_{x_0}\,$. Both jets $\,j_{x_0}^3g\,$,
$\,j_{x_0}^3\bar{g}\,$ cannot be equivalent (because the gradient
of the curvature at $\,x_0\,$ equals zero for the first metric,
whereas it is non-zero for the other one), but its differential
invariants coincide: $\,K_g(x_0)\,=\,K_{\bar{g}}(x_0)\,$ and
$\,|\mbox{grad}_{x_0}K_g|^2\,=\,|\mbox{grad}_{x_0}K_{\bar{g}}|^2\,=\,0\,$.

\smallskip

\begin{center}
\textbf{Case $\,r=4\,$.}
\end{center}

A set of generators for differential invariants of order 4\, is
given by the following five functions:
\[
p_1(j_{x_0}^4g)\,=\,K_g(x_0)\,,
\]
\[
p_2(j_{x_0}^4g)\,=\,|\mathrm{grad}_{x_0}K_g|^2\,,
\]
\[
p_3(j_{x_0}^4g)\,=\,\mathrm{trace}\,(\mathrm{Hess}_{x_0}K_g)\,,
\]
\[
p_4(j_{x_0}^4g)\,=\,\mathrm{det}\,(\mathrm{Hess}_{x_0}K_g)\,,
\]
\[
p_5(j_{x_0}^4g)\,=\,\mathrm{Hess}_{x_0}K_g(\mathrm{grad}_{x_0}K_g\,,\mathrm{grad}_{x_0}K_g)\,,
\]
where $\,\mbox{Hess}_{x_0}K_g\,:=\,(\nabla \mbox{d}K_g)_{x_0}\,$
stands for the hessian of the curvature function at $\,x_0\,$.

These above functions satisfy the following inequalities:
$$
p_2\geq 0\,\,\,\,,\,\,\,\, p_3^2-4p_4\geq 0\,\,\,\,,\,\,\,\,
(2p_5-p_2p_3)^2\leq p_2^2(p_3^2-4p_4)\,.
$$
To say it in other words, these five differential invariants
define an isomorphism of differentiable spaces
$$
\xymatrix{(p_1,\ldots,p_5):\mathbb{M}_2^4 \ar@{=}[r] & Y \subset
\mathbb{R}^5}
$$
$Y\,$ being the closed subset in $\,\mathbb{R}^5\,$ determined by
the inequalities
$$
x_2\geq 0\,\,\,\,,\,\,\,\, x_3^2-4x_4\geq 0\,\,\,\,,\,\,\,\,
(2x_5-x_2x_3)^2\leq x_2^2(x_3^2-4x_4)\,.
$$

In this case, the moduli space $\,\mathbb{M}_2^4\,$ has the
following four strata:

- The generic stratum of all classes of jets $\,j_{x_0}^4g\,$
verifying that $\,\mbox{grad}_{x_0}K_g\,$ is not an eigenvector of
$\,\mbox{Hess}_{x_0}K_g\,$ (therefore, the eigenvalues of
$\,\mbox{Hess}_{x_0}K_g\,$ are different). The type of this
stratum (group of automorphisms of its jets) is
$\,[K_1=\{\mbox{Id}\}]\,$.

- The stratum of those classes of jet metrics $\,j_{x_0}^4g\,$
verifying that $\,\mbox{grad}_{x_0}K_g\,$ is a non-zero
eigenvector of $\,\mbox{Hess}_{x_0}K_g\,$. Its type is
$\,[D_1]\,$: the group of automorphisms of each jet metric is
generated by the reflection across the vector
$\,\mbox{grad}_{x_0}K_g\,$.

- The stratum composed of those classes of jet metrics
$\,j_{x_0}^4g\,$ with $\,\mbox{grad}_{x_0}K_g\,=\,0\,$ and
verifying that the eigenvectors of $\,\mbox{Hess}_{x_0}K_g\,$ are
different. The type of this stratum is $\,[D_2]\,$: the group of
automorphisms of each jet metric is generated by the reflections
across either eigenvector of $\,\mbox{Hess}_{x_0}K_g\,$.

- The stratum of all classes of jets $\,j_{x_0}^4g\,$ with
$\,\mbox{grad}_{x_0}K_g\,=\,0\,$ and verifying that the
eigenvectors of $\,\mbox{Hess}_{x_0}K_g\,$ are both equal. The
type of the stratum is $\,[O(2)]\,$.

\bigskip

\section{Appendix A: On the notion of differential invariant of
metrics}\label{appendix}

The aim of this Appendix A is to discuss the notion of
differential invariant and to back up the Definition
\ref{definiciondeinvariante} given in Section
\ref{invariantesdiferenciales}.

The notion of differential invariant must be understood as a
particular case of the concept of regular and natural operator
between natural bundles (see \cite{kms} for an exposition of the
theory of natural bundles). What follows is an adaptation of this
point of view, getting around, though, the concept of natural
bundle.

Let $\,X\,$ be an $\,n-$dimensional smooth manifold. Let
$\,M\rightarrow X\,$ be the bundle of semi-Riemannian metrics of a
fixed signature $\,(p,q)\,$ and let $\,\mathcal{M}_X\,$ denote its
sheaf of smooth sections.

Loosely speaking, the concept of differential invariant refers to
a function ``intrinsically, locally and smoothly constructed from
a metric''. Rigorously, as it is a \textit{local} construction, a
differential invariant is a morphism of sheaves:
$$
f: \mathcal{M}_X \longrightarrow \mathcal{C}_X^\infty\,,
$$
where $\,\mathcal{C}_X^\infty\,$ stands for the sheaf of smooth
functions on $\,X\,$.

The intuition of ``intrinsic and smooth construction'' can be
encoded by saying that the morphism $\,f\,$ also satisfies the
following two properties:

\smallskip

1.- \textbf{Regularity:} If $\,\{g_s \}_{s\in S}\,$ is a family of
metrics depending smoothly on certain parameters, the family of
functions $\,\{f(g_s)\}_{s\in S}\,$ also depends smoothly on those
parameters.

To be exact, let $\,S\,$ be a smooth manifold (the space of
parameters) and let $\,U\subseteq X \times S\,$ be an open set.
For each $\,s\in S\,$, consider the open set in $\,X\,$ defined as
$\,U_s := \{ x\in X\,:\,(x,s) \in U\}\,$. A family of metrics
$\,\{g_s \in \mathcal{M}(U_s)\}_{s\in S}\,$ is said to be
\textit{smooth} if the fibre map $\,U \rightarrow S^2T^*X\,$,
$\,(x,s) \mapsto (g_s)_x\,$, is smooth. In the same way, a family
of functions $\,\{f_s \in \mathcal{C}^\infty(U_s)\}_{s\in S}\,$ is
said to be smooth if the function $\,U \rightarrow \mathbb{R}\,$,
$\,(x,s) \mapsto (f_s)(x)\,$, is smooth.

In these terms, the regularity condition expresses that for each
smooth  manifold $\,S\,$, each open set $\,U \subseteq X \times
S\,$ and each smooth family of metrics $\,\{g_s \in
\mathcal{M}(U_s) \}_{s\in S}\,$, the family of functions
$\,\{f(g_s) \in \mathcal{C}^\infty(U_s)\}_{s\in S}\,$ is smooth.

\smallskip

2.- \textbf{Naturalness:} The morphism of sheaves $\,f\,$ is
equivariant with respect to the action of local diffeomorphisms of
$\,X\,$.

That is, for each diffeomorphism $\,\tau: U \rightarrow V\,$
between open sets of $\,X\,$ and for each metric $\,g\,$ on
$\,V\,$, the following condition must be satisfied:
$$
f(\tau^* g)\,=\,\tau^*(f(g))\,.
$$

\smallskip

Taking into account the previous comments, the suitability of the
following definition is now clear:

\begin{defn} A \textbf{differential invariant} associated to
semi-Riemannian metrics (of the fixed signature) is a regular and
natural morphism of sheaves $\,f:\mathcal{M}_X \rightarrow
\mathcal{C}_X^\infty\,$.
\end{defn}

Note that this definition of differential invariant seems to be
far too general, since a differential invariant $\,f(g)\,$ is not
assumed \textit{a priori} to be constructed from the coefficients
of the metric $\,g\,$ and their subsequent partial derivatives. As
we are going to show below, this question is clarified by a
beautiful result by J. Slov\'{a}k.

For every integer $\,r\geq 0\,$, we denote by $\,J^rM\rightarrow
X\,$ the fiber bundle of $\,r-$jets of semi-Riemannian metrics on
$\,X\,$ (of the prefixed signature). The fiber bundle
$\,J^{\infty}M\rightarrow X\,$ of $\,\infty-$jets of
semi-Riemannian metrics is not a smooth manifold, but it can be
endowed with the structure of a ringed space as follows. On
$\,J^{\infty}M\rightarrow X\,$ we consider the inverse limit
topology: $\,J^{\infty}M\,=\,{\displaystyle \lim_{\leftarrow}}\,
J^rM\,$; a function on an open set $\,U\subseteq J^{\infty}M\,$ is
said to be differentiable if it is locally the composition of one
of the natural projections $\,U\subseteq J^{\infty}M \rightarrow
J^rM\,$ with a smooth function on $\,J^rM\,$. This way,
$\,J^{\infty}M\,$ is a ringed space, with its sheaf of
differentiable functions.

In a similar manner, the structure of a ringed space is defined
for the fiber of the bundle $\,J^{\infty}M \rightarrow X\,$ over a
given point $\,x_0\in X\,$:
$\,J_{x_0}^{\infty}M\,=\,{\displaystyle
\lim_{\leftarrow}}\,J_{x_0}^rM\,$.

\smallskip

\begin{thm} {\bf (Slov\'{a}k)} There exists the following
bijective correspondence:
$$
\xymatrix{\{\mathrm{differentiable\,\,functions\,\,}
\tilde{f}:J^{\infty}M
\rightarrow \mathbb{R}\} \ar@{=}[d] & \tilde{f} \ar@{|->}[d]\\
\{\mathrm{regular\,\,morphisms\,\,of\,\,sheaves\,\,}
f:\mathcal{M}_X \rightarrow \mathcal{C}_X^{\infty}\} & f}
$$
with $\,f(g)(x):=\tilde{f}(j_x^{\infty}g)\,$.
\end{thm}

\medskip

The result by Slov\'{a}k \cite{slo} refers, with a bit more of
generality, to regular morphisms between sheaves of sections of
fiber bundles.

If a regular morphism $\,\mathcal{M}_X \rightarrow
\mathcal{C}_X^{\infty}\,$ is, furthermore, natural (that is, a
differential invariant), then the corresponding smooth function
$\,\tilde{f}:J^{\infty}M\rightarrow \mathbb{R}\,$ is determined by
its restriction to the fiber $\,J_{x_0}^{\infty}M\,$ of an
arbitrary point $\,x_0\in X\,$. This assertion can be expressed
more precisely in the following way.

\smallskip

\begin{cor}\label{primercorolario} Fixed a point $\,x_0\in X\,$, the set of differential
invariants $\,f:\mathcal{M}_X \rightarrow
\mathcal{C}_X^{\infty}\,$ is in bijection with the set of
differentiable $\,\mathrm{Diff}_{x_0}-$invariant functions
$\,\tilde{f}:J_{x_0}^{\infty}M\rightarrow \mathbb{R}\,$.
\end{cor}

\medskip

\begin{defn}\label{definicionsegunda} A differential invariant
$\,f:\mathcal{M}_X\rightarrow \mathcal{C}_X^{\infty}\,$ is said to
be \textbf{of order $\,\,\leq\,\, r\,$} if the corresponding
differentiable function $\,\tilde{f}:J^{\infty}M\rightarrow
\mathbb{R}\,$ factors through the projection
$\,J^{\infty}M\rightarrow J^rM\,$.
\end{defn}

\medskip

Reformulating Corollary \ref{primercorolario} for invariants of
order $\,r\,$, we obtain that Definition \ref{definicionsegunda}
coincides with that originally given in Section
\ref{invariantesdiferenciales} (Definition
\ref{definiciondeinvariante}):

\smallskip

\begin{cor} Fixed a point $\,x_0\in X\,$, the set of all differential
invariants
$$\,f:\mathcal{M}_X \rightarrow
\mathcal{C}_X^{\infty}\,$$ of order $\,\leq r\,$ is in bijection
with the set of all smooth $\,\mathrm{Diff}_{x_0}-$invariant
functions
$$\,\tilde{f}:J_{x_0}^rM \rightarrow \mathbb{R}\,.$$
\end{cor}

\bigskip

\section{Appendix B: Classification of $\infty-$jets of
metrics}\label{appendix2}

In Section \ref{structure} we have seen that differential
invariants of order $\,\leq r\,$ classify $\,r-$jets of Riemannian
metrics at a point (Theorem \ref{clasifica}). We are now going to
generalize this result for infinite-order jets.

In the proof of next lemma we will use the following well-known
fact (\cite{bourbaki}, Chap. IX, \S\,9, Lemma 6):

Let $\,G\,$ be a compact Lie group. Every decreasing sequence of
closed subgroups $\,H_1\supseteq H_2\supseteq H_3\supseteq
\cdots\,$ stabilizes, that is, there exists an integer $\,s\,$
such that $\,H_s=H_{s+1}=H_{s+2}=\cdots$

\smallskip

\begin{lem}\label{limiteproyectivo} Let $\,G\,$ a compact Lie group and let
$$
\cdots \longrightarrow X_{r+1} \longrightarrow X_r \longrightarrow
\cdots \longrightarrow X_1
$$
be an inverse system of smooth $\,G-$equivariant maps between
smooth manifolds endowed with a smooth action of $\,G\,$. There
exists an isomorphism of ringed spaces:
$$
\begin{CD}
({\displaystyle \lim_{\leftarrow}}\,X_r)/G @= {\displaystyle \lim_{\leftarrow}}\,(X_r/G)\\
[(\ldots,x_2,x_1)] & \longmapsto & (\ldots,[x_2],[x_1])\,.
\end{CD}
$$
\end{lem}

\medskip

\begin{proof} Because of the universal quotient property,
compositions of morphisms
$$
\begin{array}{ccccc}
{\displaystyle \lim_{\leftarrow}}\,X_r & \longrightarrow & X_r & \longrightarrow & X_r/G\\
(\ldots,x_2,x_1) & \longmapsto & x_r & \longmapsto & [x_r]
\end{array}
$$
induce morphisms of ringed spaces
$$
\begin{array}{ccc}
({\displaystyle \lim_{\leftarrow}}\,X_r)/G &
\longrightarrow & (X_r/G)\\
{}[(\ldots,x_2,x_1)] & \longmapsto & [x_r]\,,
\end{array}
$$
which, for their part, because of the universal inverse limit
property, define a morphism of ringed spaces
$$
\begin{array}{ccc}
({\displaystyle \lim_{\leftarrow}}\,X_r)/G &
\stackrel{\varphi}{\longrightarrow} & {\displaystyle \lim_{\leftarrow}}\,(X_r/G)\\
{}[(\ldots,x_2,x_1)] & \longmapsto & (\ldots,[x_2],[x_1])\,.
\end{array}
$$

It is easy to check that this morphism is surjective. Let us see
that it is also injective.

First note that, given a point $\,(\ldots,x_2,x_1)\in
{\displaystyle \lim_{\leftarrow}}\,X_r$, we can get the decreasing
sequence $\,H_{x_1}\supseteq H_{x_2}\supseteq H_{x_3}\supseteq
\cdots\,$ of closed subgroups of $\,G\,$, where $\,H_{x_k}\,$
stands for the stabilizer subgroup of $\,x_k\,$. This chain
stabilizes, since $\,G\,$ is compact, so for a certain $\,s\,$ it
is verified $\,H_{x_s}=H_{x_{s+1}}=H_{x_{s+2}}=\cdots$

Let now $\,[(\ldots,x_2,x_1)]\,$ and $\,[(\ldots,x_2',x_1')]\,$ be
two points in $\,({\displaystyle \lim_{\leftarrow}}\,X_r)/G\,$
having the same image through $\,\varphi\,$, i.e.
$\,[x_k]=[x_k']\,$, for each $\,k\geq 0\,$. Write $\,x_s'=g\cdot
x_s\,$ for some $\,g\in G\,$. As the morphisms $\,X_s\rightarrow
X_k\,$ (with $\,s\geq k\,$) are $\,G-$equivariant, it is verified
that $\,x_k'=g\cdot x_k\,$ for every $\,k\leq s\,$.

Let us show that the same happens when $\,k>s\,$. As
$\,[x_k]=[x_k']\,$, we have $\,x_k'=g_k\cdot x_k\,$ for a certain
$\,g_k\in G\,$; applying that $\,X_k\rightarrow X_s\,$ is
equivariant yields $\,x_s'=g_k\cdot x_s\,$, and then (comparing
with $\,x_s'=g\cdot x_s\,$) $\,g^{-1}g_k\in H_{x_s}\,$; since
$\,H_{x_s}=H_{x_k}\,$, it follows that $\,g^{-1}g_k\in H_{x_k}\,$,
and hence the condition $\,x_k'=g_k\cdot x_k\,$ is equivalent to
$\,x_k'=g\cdot x_k\,$. In conclusion, $\,x_k'=g\cdot x_k\,$ for
every $\,k>0\,$, and therefore $\,[(\ldots,x_2,x_1)]\,$ and
$\,[(\ldots,x_2',x_1')]\,$ are the same point in
$\,({\displaystyle \lim_{\leftarrow}}\,X_r)/G\,$.

Once we have proved that $\,\varphi\,$ is bijective, it is routine
to check that $\,\varphi\,$ is an isomorphism of ringed spaces.
\end{proof}

\smallskip

\begin{defn} Let $\,x_0\in X\,$ and let
$$
J_{x_0}^{\infty}M:={\displaystyle \lim_{\leftarrow}}\,J_{x_0}^rM
$$
be the ringed space of $\infty-$jets of Riemannian metrics at
$\,x_0\,$ on $\,X\,$. The quotient ringed space
$$
\mathbb{M}_n^{\infty}:=J_{x_0}^{\infty}M/\mathrm{Diff}_{x_0}
$$
is called \textbf{moduli space} of $\,\infty-$jets of Riemannian
metrics in dimension $\,n\,$.
\end{defn}

\medskip

In the same fashion as for finite-order jets, the moduli space
$\,\mathbb{M}_n^{\infty}\,$ depends neither on the choice of the
point $\,x_0\,$ nor on that of the $\,n-$dimensional manifold
$\,X\,$.

For every integer $\,r>0\,$ we have an evident morphism of ringed
spaces
$$
\begin{array}{ccc}
\mathbb{M}_n^{\infty} & \longrightarrow & \mathbb{M}_n^r\\
{}[j_{x_0}^{\infty}g] & \longmapsto & [j_{x_0}^rg]\,,
\end{array}
$$
and these morphisms allow us to define another morphism of ringed
spaces:
$$
\begin{array}{ccc}
\mathbb{M}_n^{\infty} & \longrightarrow & {\displaystyle
\lim_{\leftarrow}}\, \mathbb{M}_n^r\\
{}[j_{x_0}^{\infty}g] & \longmapsto &
(\ldots,[j_{x_0}^rg],\ldots),.
\end{array}
$$

\smallskip

\begin{thm}\label{isomoranilladosinfinito} There exists an isomorphism of ringed spaces
$$
\begin{CD}
\mathbb{M}_n^{\infty} @= {\displaystyle
\lim_{\leftarrow}}\, \mathbb{M}_n^r\\
[j_{x_0}^{\infty}g]  & \longmapsto &
(\ldots,[j_{x_0}^rg],\ldots)\,.
\end{CD}
$$
\end{thm}

\medskip

\begin{proof} Fix a local coordinate system $\,(z_1,\ldots,z_n)\,$
centered at $\,x_0\,$. With the same notations as in Section
\ref{keylemma}, let us define
$$
\mathcal{N}^{\infty}:={\displaystyle
\lim_{\leftarrow}}\,\mathcal{N}^r\,.
$$
In other words, $\,\mathcal{N}^{\infty}\,$ is the subspace of
$\,J_{x_0}^{\infty}M\,$ formed by all those $\,\infty-$jets at
$\,x_0\,$ of Riemannian metrics having $\,(z_1,\ldots,z_n)\,$ as a
normal coordinate system. All lemmas in Section \ref{keylemma},
with their corresponding proofs, remain valid when substituting
the integer $\,\infty\,$ for $\,r\,$. In particular, our
Fundamental Lemma \ref{fundamental}, when $\,r=\infty\,$, gives us
the desired isomorphism of ringed spaces:
\[
\mathbb{M}_n^{\infty}\,=\,\left({\displaystyle \prod_{k\geq
2}}N_k\right)/O(n)\,=\,\left({\displaystyle
\lim_{\leftarrow}}\,(N_2\times\cdots\times N_r)\right)/O(n)
\]
(by Lemma \ref{limiteproyectivo})
\[
=\,{\displaystyle \lim_{\leftarrow}}\left((N_2\times\cdots\times
N_r)/O(n)\right)\,=\, {\displaystyle
\lim_{\leftarrow}}\,\mathbb{M}_n^r\,.
\]
\end{proof}

\smallskip

\begin{cor} Differential invariants of finite order classify
$\,\infty-$jets of Riemannian metrics: Two jet metrics
$\,j_{x_0}^{\infty}g\,$ and $\,j_{x_0}^{\infty}\bar{g}\,$ are
equivalent if and only if for each finite-order differential
invariant $\,h\,$ it is satisfied $\,h(g)(x_0)=h(\bar{g})(x_0)\,$.
\end{cor}

\medskip

\begin{proof} According to Theorem \ref{isomoranilladosinfinito},
we get:
$$
j_{x_0}^{\infty}g \equiv j_{x_0}^{\infty}\bar{g}\,\,\,
\Longleftrightarrow \,\,\, j_{x_0}^rg\equiv
j_{x_0}^r\bar{g}\,\,\,,\,\,\,\forall r \geq 0\,.
$$

To complete our proof, it is sufficient to use the fact that
differential invariants of order $\,\leq r\,$ classify $\,r-$jet
metrics (Theorem \ref{clasifica}).
\end{proof}

\end{document}